\documentclass[12pt]{article}

\usepackage[latin1]{inputenc}

\usepackage{pgf,tikz}
\usetikzlibrary{arrows}
\usetikzlibrary{patterns}

\usepackage{multirow} 
\usepackage{subfigure}

\usepackage{amssymb}

\usepackage{amsmath,epsfig}
\def\bar{\textrm{bar}}

\newtheorem{Theorem}{Theorem}
\newtheorem{Corollary}[Theorem]{Corollary}

\newtheorem{Definition}[Theorem]{Definition}
\newtheorem{Lemma}[Theorem]{Lemma}

\numberwithin{equation}{section}

\newcommand{\N}{{\mathbb N}}

\newcommand{\C}{{\mathbb C}}

\newcommand{\beq}{\begin{equation}}
\newcommand{\eeq}{\end{equation}}

\def\qed{\hfill$\Box$\\ \medskip}

\usepackage{multirow}

\def\mybox{\hbox to 12.0pt}

\def\mybigbox{\hbox to 35.0pt}
\def\myverybigbox{\hbox to 60.0pt}

\def\ds{\displaystyle}    
\def\ov{\overline}
\def\sc{\scriptstyle}

\def\wgt{{\rm wgt}}

\def\ch{{\rm ch}\,}
\def\w{{\bf w}}
\def\x{{\bf x}}

\def\z{{\bf z}}

\def\a{{\bf a}}

\def\0{{\bf 0}}

\def\bfhalf{{\bf \frac12}}

\def\ovx{{\bf\ov{\x}}}

\def\g{g}
\def\gl{gl}
\def\oo{oo}
\def\sp{sp}
\def\eo{eo}
\def\spd{spd}
\def\eod{eod}

\def\red{\textcolor{red}}
\def\blue{\textcolor{blue}}
\def\magenta{\textcolor{magenta}}
\def\cyan{\textcolor{cyan}}

\newcommand{\YT}[4]{
\vcenter{\hbox{
\begin{tikzpicture}[x={(0in,-#1)},y={(#2,0in)}] 
\foreach \k [count=\i] in {#4} {
 \foreach \e [count=\j] in \k {
  \draw (\i,\j) rectangle +(-1,-1);
  \draw (\i-0.5,\j-0.5) node {$#3\e$};
 }
}
\end{tikzpicture}
}}
}

\newcommand{\wideYT}[4]{
\vcenter{\hbox{
\begin{tikzpicture}[x={(0in,-#1)},y={(#2,0in)}] 
\foreach \k [count=\i] in {#4} {
 \foreach \e [count=\j] in \k {
  \draw (\i,\j) rectangle +(-1,-1);
  \draw (\i-0.5,\j-0.5) node {$#3\e$};
 }
}
\end{tikzpicture}
}}
}

\title{Factorial characters of the classical Lie groups and their combinatorial realisations}

\author{
Ang\`ele M. Hamel
\thanks{ Department of Physics and Computer Science,
Wilfrid Laurier University,
Waterloo, Ontario, N2L 3C5, Canada ({\tt ahamel@wlu.ca})}
\and 
Ronald C. King\thanks{
Mathematical Sciences, University of Southampton, 
Southampton SO17 1BJ, England ({\tt R.C.King@soton.ac.uk})}}

\begin{document}

\maketitle

\begin{abstract}
Just as the definition of factorial Schur functions as a ratio of determinants 
allows one to show that they satisfy a Jacobi-Trudi-type identity
and have an explicit combinatorial realisation in terms of semistandard tableaux, 
so we offer here definitions of factorial irreducible characters of the classical Lie groups as ratios of determinants 
that share these two features. These factorial characters are each specified by a partition, 
$\lambda=(\lambda_1,\lambda_2,\ldots,\lambda_n)$, and in each case a flagged Jacobi-Trudi identity 
is derived that expresses the factorial character as a determinant 
of corresponding factorial characters specified by one-part partitions, $(m)$,
for which we supply generating functions. 
These identities are established by manipulating determinants through the use 
of certain recurrence relations derived from 
these generating functions. The transitions to
combinatorial realisations of the factorial characters in terms of tableaux are then established by means of non-intersecting lattice path models.
The results apply to $gl(n)$, $so(2n+1)$, $sp(2n)$ and $o(2n)$, and are extended to the case of $so(2n)$
by making use of newly defined factorial difference characters.
\end{abstract}

\section{Introduction}
Let $n\in\N$ be fixed. Let $\x=(x_1,x_2,\ldots,x_n)$ and $\ovx=(\ov{x}_1,\ov{x}_2,\ldots,\ov{x}_n)$
with $\ov{x}_i=x_i^{-1}$ for $i=1,2,\ldots,n$, and let $\lambda=(\lambda_1,\lambda_2,\ldots,\lambda_n)$ be a partition 
of length $\ell(\lambda)\leq n$. Then each of the classical groups $G=GL(n,\C)$, $SO(2n+1,\C)$, $Sp(2n,\C)$,  
and $O(2n,\C)$ possesses a finite dimensional irreducible representation $V_G^\lambda$ of highest weight $\lambda$ whose
character may be denoted by $\ch V_G^\lambda(\z)$, where $\z$ is a suitable parametrisation of the eigenvalues
of the group elements of $G$, namely $\x$, $(\x,\ovx,1)$, $(\x,\ovx)$ and $(\x,\ovx)$, respectively. 
The irreducible representation of $O(2n,\C)$ of highest weight $\lambda$ remains irrreducible on restriction 
to $SO(2n,\C)$ if $\ell(\lambda)<n$, but reduces to a direct sum of a pair of irreducible representation of $SO(2n,\C)$
of highest weights $\lambda_+$ and $\lambda_-$ if $\ell(\lambda)=n$, where $\lambda_+=(\lambda_1,\ldots,\lambda_{n-1},\lambda_n)$
and $\lambda_-=(\lambda_1,\ldots,\lambda_{n-1},-\lambda_n)$.

For the sake of typographical simplicity we adopt the following notation for the above characters:
\begin{align}\label{eqn-char-sfn}
&\ch V_{GL(n,\C)}^\lambda(\x) = gl_\lambda(\x)\,; \qquad \ch V_{SO(2n+1,\C)}^\lambda(\x,\ovx,1) =  so_\lambda(\x,\ovx,1)\,; \cr
&\ch V_{Sp(2n,\C)}^\lambda(\x,\ovx) = sp_\lambda(\x,\ovx)\,; \qquad \ch V_{O(2n,\C)}^\lambda(\x,\ovx) = o_\lambda(\x,\ovx)\,.
\end{align}
with
\begin{equation}
o_\lambda(\x,\ovx)= \begin{cases} so_\lambda(\x,\ovx) &\hbox{if $\lambda_n=0$};\cr
   so_{\lambda_{+}}(\x,\ovx)+so_{\lambda_{-}}(\x,\ovx)&\hbox{if $\lambda_n>0$}.\cr \end{cases} \nonumber
\end{equation}
There also exist difference characters of $SO(2n,\C)$ defined by
\begin{equation}
o'_\lambda(\x,\ovx)= \begin{cases} 0&\hbox{if $\lambda_n=0$};\cr
   so_{\lambda_{+}}(\x,\ovx)-so_{\lambda_{-}}(\x,\ovx)&\hbox{if $\lambda_n>0$}.\cr \end{cases} \nonumber
\end{equation}
The above notation reflects the fact that, with a slightly different interpretation of the various parameters $x_i$ and $\ov{x}_i$
in $\z$ as formal exponentials $e^{\epsilon_i}$ and $e^{-\epsilon_i}$ of Euclidean basis vectors $\epsilon_i$ in the 
weight space of the Lie algebras $\g=gl(n)$, $sp(2n)$, $so(2n+1)$ and $o(2n)$, 
each $\g_\lambda(\z)$ is a character of the Lie algebra $\g$. As a convenient mnenomic we use superscripts $\gl$, $\sp$, $\oo$ and $\eo$ 
to distinguish our general linear, symplectic, odd orthogonal and even orthogonal cases.

By means of rather simple modifications of Weyl's formula for each  character $\g_\lambda(\z)$ we offer
definitions of factorial characters $\g_\lambda(\z|\a)$ involving an infinite sequence $\a=(a_1,a_2,\ldots)$ of
completely arbitrary factorial parameters. That these definitions are appropriate is then attested to by showing 
that each $\g_\lambda(\z|\a)$ satisfies a flagged Jacobi-Trudi identity of the form
$\g_\lambda(\z|\a)=|\, \g_{(\lambda_j-j+i)}(\z^{(i)}|\a)\,|$ where the various $\z^{(i)}$ for $i=1,2,\ldots,n$ 
constitute a flag. We then provide a lattice path interpretation of each one part partition factorial character 
$\g_{(m)}(\z^{(i)}|\a)$ and by a careful specification of the underlying lattice for each $\g$ we are able to express
$\g_\lambda(\z|\a)$ itself in terms of suitable edge weighted $n$-tuples of non-intersecting lattice paths. This then allows 
us to express $\g_\lambda(\z|\a)$ as a sum over tableaux $T\in{\cal T}_\lambda^\g$ of shape $\lambda$, 
whose entries carry weights of the form $z_k+a_\ell$, or in some instances $z_k-a_\ell$, for some $k$ and $\ell$. 
The $\g$-dependent sets of tableaux ${\cal T}_\lambda^\g$ coincide with those known classically, so that all it needs to recover 
classical characters is to set $\a=\0=(0,0,\ldots)$.

Our starting point is therefore 
Weyl's character formula for $\g_\lambda(\z)$ for the appropriate Lie algebra $\g$ and the appropriate $\z$ as specified above.
It takes the form of ratio of determinants, or, in the even orthogonal case, 
as a sum of two such ratios with the same denominator. Explicit formulae, that may be taken for example from~\cite{Lit50} or~\cite{FH04},
are provided in Section~\ref{sec-fchar} for the cases $gl(n)$, $sp(2n)$, $so(2n+1)$ and $o(2n)$, while we defer those appropriate
to irreducble and difference characters of $so(2n)$ to Section~\ref{sec-so2n}. 
They are all notable for the fact that each element in each determinant 
is of the form $x^m$ or $x^{m+\epsilon} \pm \ov{x}^{m+\epsilon}$ for some $x=x_i$ with $i\in\{1,2,\ldots,n\}$
and some non-negative integer $m$, with $\epsilon=1/2$, $1$ or $0$, depending on the Lie algebra in question.

Beidenharn and Louck~\cite{BL89} introduced the notion of a factorial Schur function defined in terms of 
the Gelfand-Tsetlin patterns associated with $GL(n,\C)$.
This was further studied by Chen and Louck~\cite{CL93} before it was given its more general form by Goulden and Greene~\cite{GG94} and Macdonald~\cite{Mac92}, 
expressed this time in terms of semistandard tableaux, with Macdonald also giving an alternative definition as a ratio of alternants. It is 
this latter form that we wish to generalise to the case of the other classical Lie groups $G$.

In the $GL(n,\C)$ case the transition from the Schur function $s_\lambda(\x)$ to the factorial Schur function, $s_\lambda(\x|\a)$, 
involves an infinite sequence of factorial parameters $\a=(a_1,a_2,\ldots)$. The transition is effected
by replacing each non-negative power $x^m$ in the determinantal expression for $\gl_\lambda(\x)$ by its factorial power $(x|\a)^m$ defined by:
\begin{equation}\label{eqn-xam}
    (x|\a)^m = \begin{cases} (x+a_1)(x+a_2)\cdots(x+a_m) &\hbox{if $m>0$};\cr 1&\hbox{if $m=0$}. \cr\end{cases}
\end{equation}
In the case of the other classical Lie groups we propose that the transition from the ordinary character 
$\g_\lambda(\x)$ to the factorial character $\g_\lambda(\x|\a)$ is effected by replacing not only each $x^m$ by $(x|\a)^m$
but also each $\ov{x}^m$ by $(\ov{x}|\a)^m$ with  
\begin{equation}\label{eqn-ovxam}
    (\ov{x}|\a)^m = \begin{cases} (\ov{x}+a_1)(\ov{x}+a_2)\cdots(\ov{x}+a_m) &\hbox{if $m>0$};\cr 1&\hbox{if $m=0$}. \cr\end{cases}
\end{equation}
This results in our proposed formulae for factorial characters $\g_\lambda(\x|\a)$ given in Definition~\ref{def-fchar}.

Following Macdonald~\cite{Mac92}, it is useful to introduce a shift operator $\tau$ defined in such a way that
\begin{equation}\label{eqn-tau}
    \tau^r \a = (a_{r+1},a_{r+2},\ldots) \quad\hbox{for any integer $r$ and any $\a=(a_1,a_2,\ldots)$}\,.
\end{equation}	
This necessarily extends the sequence $\a=(a_1,a_2,\ldots)$ of factorial parameters to a doubly
infinite sequence $(\ldots,a_{-2},a_{-1},a_0,a_1,a_2,\ldots)$. 

That the definitions are appropriate depends to what extent the properties of factorial characters are analogous to
those of factorial Schur functions. We have in mind things like deriving for each of our factorial characters a factorial 
Jacobi-Trudi-type identity and a combinatorial interpretation in terms of tableau, both of which were established in the case of factorial Schur 
functions by Macdonald~\cite{Mac92,Mac95}, and perhaps more ambitiously, the derivation of Tokuyama type identities~\cite{Tok88} 
as recently derived in the factorial Schur function case by Bump, McNamara and Nakasuji~\cite{BMN14}, with an alternative derivation 
appearing in~\cite{HK15b}. This latter step necessarily requires generalisations of the classical Schur $Q$-functions~\cite{Sch11, Ste89} 
first from the general linear case to that of the other classical Lie groups, and then to factorial versions of these. 
This second part of our programme is deferred to a separate publication~\cite{HK17}.

The key to accomplishing the first part of our programme is the identification of appropriate analogues of the 
complete homogeneous symmetric functions $h_m(\x)$ and of their factorial counterparts $h_m(\x|\a)$. 
This is done for each group in Section~\ref{sec-fJT} by means of generating functions. 
A sequence of Lemmas in Section~\ref{sec-fJT} then leads to our first main result,
Theorem~\ref{thm-fJT}, which provides flagged Jacobi-Trudi identities for the factorial characters of each Lie algebra $\g$ in the form:
\begin{equation}\label{eqn-g-fJT}
  \g_\lambda(\z|\a) = \big|\, h^\g_{\lambda_j-j+i}(\z^{(i)}|\a)\,\big|\,,
\end{equation} 
where $\z$ and the precise form of $\z^{(i)}$ vary from case to case. These results amount to a factorial generalisation
of the flagged Jacobi-Trudi formulae for characters of the classical Lie groups first provided by Okada~\cite{Oka??}. Inspired in part by
preliminary presentations~\cite{HK16a,HK16b}
of the results given here, Okada~\cite{Oka16a} has also obtained by slightly 
different means not only the same flagged Jacobi-Trudi identities, but also flagged dual Jacobi-Trudi and unflagged Giambelli identities 
for all the factorial characters $\g_\lambda(\z|\a)$.

In Section~\ref{sec-hma} the generating function definitions of $h^{\g}_m(\x|\a)$ are used to establish recurrence relations
that enable each of them to be expressed as sum of products of factors of the form $x_k+a_\ell$ and $\ov{x}_k+a_\ell$. This allows us 
to interpret each term contributing to $h^{\g}_m(\x|\a)$ as a suitable weighted sequence of edges constituting a path in a $\g$-dependent
lattice. Each determinant of (\ref{eqn-g-fJT}) is then evaluated in Section~\ref{sec-tableaux} as a signed sum of contributions from
$n$-tuples of such lattice paths, that can be reduced to a sum over $n$-tuples of non-intersecting lattice paths by means of the 
Lindstr\"om-Gessel-Viennot Theorem~\cite{Lin73,GV85,GV89}. Again for each $\g$, a bijective correspondence between such $n$-tuples and
tableaux leads to our second main result, namely Theorem~\ref{thm-Twgts}, which offers a combinatoral realisation
of factorial characters of the form
\begin{equation}\label{eqn-g-Twgts}
   \g_\lambda(\z|\a) = \sum_{T\in{\cal T}^{\g}_\lambda}\ 2^{\zeta(T)}\prod_{(i,j)\in F^\lambda} \wgt(T_{ij})\,.
\end{equation}
Here the sum is over a set ${\cal T}^{\g}$ of tableaux $T$ consisting of arrays of entries in a Young diagram $F^\lambda$ of shape $\lambda$.
The entries, taken from alphabets appropriate to each $\g$, are subject to a variety of conditions such as the semistandardness condition that
applies in the $\gl(n)$ case. The sets ${\cal T}^{\g}$ coincide with familiar sets of classical group tableaux~\cite{Lit50,KElS83,Sun90}.
Here $\wgt(T_{ij})$ is the factorial weight of the entry $T_{ij}$ at position $(i,j)$ in $F^\lambda$ and generally takes the form
of $(x_k+a_\ell)$ and $(\ov{x}_k+a_\ell)$ for some $k$ and $\ell$. The factor $2^{\zeta(T)}$ is peculiar to the $o(2n)$ case.

In Section~\ref{sec-so2n} the results are extended to the case of $so(2n)$ for which complications arise as a result of restricting 
some irreducible representations of $O(2n,\C)$ to its subgroup $SO(2n,\C)$. These complications are dealt with by invoking the 
difference characters of $so(2n)$ introduced by Littlewood~\cite{Lit50}. The factorial version of these characters are shown
to satisfy a flagged Jacobi-Trudi identity. Proceeding exactly as in the case of $o(2n)$ factorial characters we construct non-intersecting path 
and tableaux models, enabling us to evaluate factorial difference characters. Combining the results for ordinary and difference characters
leads to our final main result Theorem~\ref{thm-Twgts-so2n}, expressing factorial irreducible characters of $so(2n)$ as sums over suitably 
weighted tableaux. We offer some brief concluding remarks in Section~\ref{sec-remarks}.

\section{Ordinary and factorial characters}
\label{sec-fchar}

The first of our characters, $gl_\lambda(\x)$, in (\ref{eqn-char-sfn}) is none other than the Schur polynomial $s_\lambda(\x)$ which has a well known 
definition as a ratio of alternants~\cite{Lit50,Mac95}. Thanks to Weyl's character formula for the corresponding Lie algebras 
each of the our characters under consideration can be expressed in a similar form~\cite{FH04}:
\begin{align}
\gl_\lambda(\x) 
&= \frac{\left|\, x_i^{\lambda_j+n-j} \,\right|}
				{\left|\, x_i^{n-j} \,\right| }\,; \label{eqn-gl-char}\\
sp_\lambda(\x,\ovx) 
&= \frac{\left|\, x_i^{\lambda_j+n-j+1} - \ov{x}_i^{\lambda_j+n-j+1} \,\right|}
				                          {\left|\, x_i^{n-j+1} - \ov{x}_i^{n-j+1} \,\right|}\,; \label{eqn-sp-char}\\
so_\lambda(\x,\ovx,1)
&=  \frac{\left|\, x_i^{\lambda_j+n-j+1/2} - \ov{x}_i^{\lambda_j+n-j+1/2} \,\right|}
				                          {\left|\, x_i^{n-j+1/2} - \ov{x}_i^{n-j+1/2} \,\right|}\,; \label{eqn-oo-char}\\																	
o_\lambda(\x,\ovx)
&=  \frac{\eta\left|\, x_i^{\lambda_j+n-j} + \ov{x}_i^{\lambda_j+n-j} \,\right| }
				                          {\frac12\left|\, x_i^{n-j} + \ov{x}_i^{n-j} \,\right|} 
																	 ~~\hbox{with}~~\eta=\begin{cases}\frac12&\hbox{if $\lambda_n=0$};\cr
																	                                  1&\hbox{$\lambda_n>0$}.
																																		    \end{cases}. \label{eqn-eo-char}															
\end{align}
where all the determinants, both here and hereafter, are $n\times n$, and we have specified in each case the 
element in the $i$th row and $j$th column.

Weyl's denominator formula~\cite{FH04} is such that
\begin{align}
\left|\, x_i^{n-j} \,\right| &= \prod_{1\leq i<j\leq n} (x_i-x_j)\,; \label{eqn-gl-denom}\\
\left|\, x_i^{n-j+1} - \ov{x}_i^{n-j+1} \,\right| &= \prod_{i=1}^n (x_i-\ov{x}_i)\ \prod_{1\leq i<j\leq n} (x_i+\ov{x}_i-x_j-\ov{x}_j)\,;  \label{eqn-sp-denom}\\
\left|\, x_i^{n-j+1/2} - \ov{x}_i^{n-j+1/2} \,\right| &= \prod_{i=1}^n (x_i^{1/2}-\ov{x}_i^{1/2})\ \prod_{1\leq i<j\leq n} (x_i+\ov{x}_i-x_j-\ov{x}_j)\,; \label{eqn-oo-denom} \\
{\sc{\frac12}}\left|\, x_i^{n-j} + \ov{x}_i^{n-j} \,\right| &= \prod_{1\leq i<j\leq n} (x_i+\ov{x}_i-x_j-\ov{x}_j)\,. \label{eqn-eo-denom}
\end{align}

As anticipated in the Introduction, we propose the following definition of factorial characters of the classical Lie algebras
based on the use of $(\x|\a)^m$ and $(\ov{x}|\a)^m$ as given in (\ref{eqn-xam}) and (\ref{eqn-ovxam}):
\begin{Definition}\label{def-fchar}
For any partition $\lambda=(\lambda_1,\lambda_2,\ldots,\lambda_n)$ of length $\ell(\lambda)\leq n$ let 
\begin{align}
gl_\lambda(\x|\a) 
&= \frac{\left|\, (x_i|\a)^{\lambda_j+n-j} \,\right|}
				{\left|\, (x_i|\a)^{n-j} \,\right| }\,; \label{eqn-gla}\\
sp_\lambda(\x,\ovx|\a) 
&= \frac{\left|\, x_i(x_i|\a)^{\lambda_j+n-j} - \ov{x}_i(\ov{x}_i|\a)^{\lambda_j+n-j} \,\right|}
				                          {\left|\, x_i(x_i|\a)^{n-j} - \ov{x}_i(\ov{x}_i|\a)^{n-j} \,\right|}\,; \label{eqn-spa}\\
so_\lambda(\x,\ovx,1|\a)
&=  \frac{\left|\, x_i^{1/2} (x_i|\a)^{\lambda_j+n-j} - \ov{x}_i^{1/2}(\ov{x}_i|\a)^{\lambda_j+n-j} \,\right|}
				                          {\left|\, x_i^{1/2}(x_i|\a)^{n-j} - \ov{x}_i^{1/2}(\ov{x}_i|\a)^{n-j} \,\right|}\,; \label{eqn-ooa}\\
o_\lambda(\x,\ovx|\a)
&=  \frac{\eta\left|\, (x_i|\a)^{\lambda_j+n-j} + (\ov{x}_i|\a)^{\lambda_j+n-j} \,\right|}
				                          {\frac12\left|\, (x_i|\a)^{n-j} + (\ov{x}_i|\a)^{n-j} \,\right|}
																	 ~~\hbox{with}~~\eta=\begin{cases}\frac12&\hbox{if $\lambda_n=0$};\cr
																	                                  1&\hbox{$\lambda_n>0$}.
																																		    \end{cases}. \label{eqn-eoa}																	
\end{align}
\end{Definition}

The first definition is that of Macdonald~\cite{Mac95} for factorial Schur functions, and the others have been drawn up as rather 
natural generalisations of this that all have the merit of reducing to the classical non-factorial characters if one sets $\a=\0=(0,0,\ldots)$.
In each case the denominators are independent of $\a$ and coincide with the Weyl denominators of (\ref{eqn-gl-denom})-(\ref{eqn-eo-denom}),
as can be seen by extending the argument of Macdonald~\cite{Mac95} and noting that $(x_i|\a)^{n-j}$ and $(\ov{x}_i|\a)^{n-j}$ are monic polynomials
of degree $n-j$ with identical coefficients of powers $x_i^k$ and $\ov{x}_i^k$ for all $k$, where the exponents $n-j$ take the values $n,\ldots,1,0$ across 
consecutive columns of each determinant. Successively subtracting suitable $\a$-dependent multiples of columns from right to left eliminates
all but those terms independent of $\a$.

\section{Preliminary lemmas and a flagged Jacobi-Trudi theorem}
\label{sec-fJT}

To establish factorial Jacobi-Trudi identites we need analogues of the complete homogeneous
symmetric functions $h_r(\x)$ that are appropriate not only to the case of the other group characters but also to 
the case of our factorial characters. Just is done classically for $h_r(\x)$ it is convenient to define these analogues 
by means of generating functions. Our notation is such that each generating function $F(\z;t)$ may be expanded
as a power series in $t$, and we denote the coefficient of $t^m$ in such an expansion by $[t^m]\ F(\z;t)$ for all
integers $m$.

\begin{Definition}\label{def-hm}
For any integer $m$ let 
\begin{align}
h^{\gl}_m(\x) &= [t^m]\ \prod_{i=1}^n \frac{1}{1-tx_i}\,;\label{eqn-gl-hm} \\
h^{\sp}_m(\x,\ovx) &= [t^m]\ \prod_{i=1}^n \frac{1}{(1-tx_i)(1-t\ov{x}_i)}\,;\label{eqn-sp-hm} \\
h^{\oo}_m(\x,\ovx,1) &= [t^m]\ (1+t)\ \prod_{i=1}^n \frac{1}{(1-tx_i)(1-t\ov{x}_i)}\,;\label{eqn-oo-hm} \\
h^{\eo}_m(\x,\ovx) &= [t^m]\ \begin{cases} \ds \left(\frac{1}{1-tx_1}+\frac{1}{1-t\ov{x}_1}-\delta_{m0}\right)&\hbox{if $n=1$}\,;\cr
                            \ds  (1-t^2)\ \prod_{i=1}^n \frac{1}{(1-tx_i)(1-t\ov{x}_i)}&\hbox{if $n>1$}\,,\cr \end{cases}
														\label{eqn-eo-hm}
\end{align}
with $\delta_{m0}=1$ if $m=0$ and $0$ otherwise.
In the case $m=0$ we have $h^{\gl}_0(\x)= h^{\sp}_0(\x,\ovx) = h^{\oo}_0(\x,\ovx,1) = h^{\eo}_0(\x,\ovx) = 1$,
while for $m<0$ we have $h^{\gl}_m(\x)= h^{\sp}_m(\x,\ovx) = h^{\oo}_m(\x,\ovx,1) = h^{\eo}_m(\x,\ovx) = 0$.
\end{Definition} 

In the factorial case we make the following definitions in terms of generating functions $F_m(\z|\a;t)$ that are truncated 
in the sense that the power $m$ of $[t^m]$ appears in an upper limit of the associated generating function.

\begin{Definition}\label{def-hma}
For any integer $m$ let 
\begin{align}
h^{\gl}_m(\x|\a) &= [t^m]\ \prod_{i=1}^n \frac{1}{1-tx_i} \ \prod_{j=1}^{n+m-1}(1+ta_j)\,;\label{eqn-gl-hma} \\
h^{\sp}_m(\x,\ovx|\a) &= [t^m]\ \prod_{i=1}^n \frac{1}{(1-tx_i)(1-t\ov{x}_i)}\ \prod_{j=1}^{n+m-1}(1+ta_j)\,;\label{eqn-sp-hma} \\
h^{\oo}_m(\x,\ovx,1|\a) &= [t^m]\ (1+t)\ \prod_{i=1}^n \frac{1}{(1-tx_i)(1-t\ov{x}_i)}\ \prod_{j=1}^{n+m-1}(1+ta_j)\,;\label{eqn-oo-hma} \\
h^{\eo}_m(\x,\ovx|\a) &= [t^m]\ \begin{cases} \ds \left(\frac{1}{1-tx_1}+\frac{1}{1-t\ov{x}_1}-\delta_{m0}\right)\ \prod_{j=1}^{m}(1+ta_j)&\hbox{if $n=1$}\,;\cr
                            \ds  (1-t^2)\ \prod_{i=1}^n \frac{1}{(1-tx_i)(1-t\ov{x}_i)}\ \prod_{j=1}^{n+m-1}(1+ta_j)&\hbox{if $n>1$}\,.\cr \end{cases}
														\label{eqn-eo-hma}
\end{align}
Then for $m=0$ we have $h^{\gl}_0(\x|\a)= h^{\sp}_0(\x,\ovx|\a) = h^{\oo}_0(\x,\ovx,1|\a) = h^{\eo}_0(\x,\ovx|\a) = 1$,
while for $m<0$ we have $h^{\gl}_m(\x|\a)= h^{\sp}_m(\x,\ovx|\a) = h^{\oo}_m(\x,\ovx,1|\a) = h^{\eo}_m(\x,\ovx|\a) = 0$.
\end{Definition} 

It might be noted here, that each of these expressions is manifestly symmetric not only 
under all permutations of the parameters $x_i$ and $\ov{x}_i$ for $i=1,2,\ldots,n$, 
including interchanges of $x_i$ and $\ov{x}_j$ for any $i,j=1,2,\ldots,n$,
but also under all permutations of the parameters $a_j$ for $j=1,2,\ldots,n+m-1$.

\bigskip

The one variable cases $\x=(x_i)$ of Definition~\ref{def-hma} allow us to rewrite
our factorial characters in the following manner:
\begin{Lemma}\label{Lem-char-hmxi}
For any partition $\lambda=(\lambda_1,\lambda_2,\ldots,\lambda_n)$ of length $\ell(\lambda)\leq n$ we have
\begin{align}
\gl_\lambda(\x|\a) 
&= \frac{\left|\, h^{\gl}_{\lambda_j+n-j}(x_i|\a) \,\right|}
				{\left|\, h^{\gl}_{n-j}(x_i|\a) \,\right| }\,; \label{eqn-gl-hh}\\						
sp_\lambda(\x,\ovx|\a) 
&= \frac{\left|\, h^{\sp}_{\lambda_j+n-j}(x_i,\ov{x}_i|\a) \,\right|}
				                          {\left|\, h^{\sp}_{n-j}(x_i,\ov{x_i}|\a) \,\right|}\,; \label{eqn-sp-hh}\\
so_\lambda(\x,\ovx,1|\a)
&=  \frac{\left|\, h^{\oo}_{\lambda_j+n-j}(x_i,\ov{x}_i,1|\a) \,\right|}
				                          {\left|\, h^{\oo}_{n-j}(x_i,\ov{x}_i|\a) \,\right|}\,; \label{eqn-oo-hh}\\
o_\lambda(\x,\ovx|\a)
&= \frac{\left|\, h^{\eo}_{\lambda_j+n-j}(x_i,\ov{x}_i|\a) \,\right|} 
				                          {\left|\, h^{\eo}_{n-j}(x_i,\ov{x}_i|\a) \,\right|}\,. \label{eqn-eo-hh}																	
\end{align}
\end{Lemma}

\noindent{\bf Proof}:
In the case of $\gl_\lambda(\x|\a)$ it suffices to note that for $m\geq 0$
\begin{align}\label{eqn-gl-hxam}
    h_m^{\gl}(x_i|\a) &= [t^m] \frac{1}{1-tx_i} \prod_{j=1}^m (1+ta_j) =  [t^m] \frac{1+ta_1}{1-tx_i} \prod_{j=2}^{m} (1+ta_j)\cr
                &= [t^m] \left(1+\frac{t(x_i+a_1)}{1-tx_i}\right) \prod_{j=2}^{m} (1+ta_j) = (x_i+a_1) [t^{m-1}] \frac{1}{1-tx_i} \prod_{j=2}^{m} (1+ta_j)\cr
								&= (x_i+a_1)(x_i+a_2)\cdots(x_i+a_m) [t^0] \frac{1}{1-tx_i} = (x_i|\a)^m\,.
\end{align}
One then just uses this identity in (\ref{eqn-gla}) with $m=\lambda_j+n-j$ and $m=n-j$
in the numerator and denominator as appropriate.

In the $o_\lambda(\x,\ovx,1|\a)$ case we have
\begin{align}\label{eqn-oo-hxma}
h_m^{\oo}(x_i,\ov{x}_i,1|\a)&=[t^m]\frac{1+t}{(1-tx_i)(1-t\ov{x}_i)} \prod_{j=1}^m (1+ta_j)\cr
	          &=[t^m] \frac{1}{x_i^{1/2}-\ov{x}_i^{1/2}}\left(\frac{x_i^{1/2}}{1-tx_i}-\frac{\ov{x}_i^{1/2}}{1-t\ov{x}_i}\right) \prod_{j=1}^m (1+ta_j)\cr
	                    &=\frac{1}{x_i^{1/2}-\ov{x}_i^{1/2}}\left( x_i^{1/2}(x_i|\a)^m-\ov{x}_i^{1/2}(\ov{x}_i|\a)^m\right)\,,
\end{align}
where use has been made of (\ref{eqn-gl-hxam}).
Then one again uses this identity in the numerator and denominator of (\ref{eqn-ooa}) with $m=\lambda_j+n-j$ and $m=n-j$ as appropriate,
and cancels the common factors of $x_i^{1/2}-\ov{x}_i^{1/2}$ that can be extracted from the 
$i$th row in both numerator and denominator for $i=1,2,\ldots,n$.

Similarly, in the $sp_\lambda(\x,\ovx|\a)$ case we have 
\begin{align}\label{eqn-sp-hxma}
h_m^{\sp}(x_i,\ov{x}_i|\a)&=[t^m]\frac{1}{(1-tx_i)(1-t\ov{x}_i)} \prod_{j=1}^m (1+ta_j) \cr
	          &=[t^m] \frac{1}{x_i-\ov{x}_i}\left(\frac{x_i}{1-tx_i}-\frac{\ov{x}_i}{1-t\ov{x}_i}\right) \prod_{j=1}^m (1+ta_j)\cr
	                    &=\frac{1}{x_i-\ov{x}_i}\left( x_i(x_i|\a)^m-\ov{x}_i(\ov{x}_i|\a)^m\right)\,.
\end{align}
The required result follows as before with the cancellation this time of the common factors $x_i-\ov{x}_i$.

Finally, the case of  $o_\lambda(\x,\ovx|\a)$ only depends on the $n=1$ case of (\ref{eqn-eo-hma}). This yields
\begin{align}\label{eqn-eo-hxma}
   h_m^{\eo}(x_i,\ov{x}_i|\a)&=[t^m]\left(\frac{1}{1-tx_i}+\frac{1}{1-t\ov{x}_i}-\delta_{m0}\right)\prod_{j=1}^m (1+ta_j)\cr
									 &=(x_i|\a)^m+(\ov{x}_i|\a)^m-\delta_{m0}\,,
\end{align}
as required to complete the proof by exploiting the cases $m=\lambda_j+n-j$ and $m=n-j$. 
The term $-\delta_{m0}$ ensures that the factors $\eta$ and $\frac12$ that 
appear in (\ref{eqn-eoa}) do not appear in (\ref{eqn-eo-hma}).
\qed

The next step is to transform each of the expressions in Lemma~\ref{Lem-char-hmxi} into flagged Jacobi-Trudi identities.
This is accomplished by means of the following Lemma:

\begin{Lemma}\label{Lem-hra-xixj}
Let $\x=(x_1,x_2,\ldots,x_n)$, $\ovx=(\ov{x}_1,\ov{x}_2,\ldots,\ov{x}_n)$ and $\a=(a_1,a_2,\ldots)$, including the case
$\a=\0=(0,0,\ldots)$. Then for $1\leq i<j\leq n$ and all integers $m$:
\begin{align}
  &h^{\gl}_m(x_i,\ldots x_{j-1}|\a)-h^{\gl}_m(x_{i+1},\ldots,x_j|\a)\cr
	&\hskip4cm=(x_i-x_j) h^{\gl}_{m-1}(x_i,\ldots,x_j|\a)\,; \label{eqn-gl-hra-xixj} \\
	&h^{\sp}_m(x_i,\ov{x}_i,\ldots,x_{j-1},\ov{x}_{j-1}|\a)-h^{\sp}_m(x_{i+1},\ov{x}_{i+1},\ldots,x_{j},\ov{x}_{j},|\a)\cr
	&\hskip4cm=(x_i+\ov{x}_i-x_j-\ov{x}_j) h^{\sp}_{m-1}(x_i,\ov{x}_i,\ldots,x_{j},\ov{x}_{j},|\a)\; \label{eqn-sp-hra-xixj} \\
	&h^{\oo}_m(x_i,\ov{x}_i,\ldots,x_{j-1},\ov{x}_{j-1},1|\a)-h^{\oo}_m(x_{i+1},\ov{x}_{i+1},\ldots,x_{j},\ov{x}_{j},1|\a)\cr
	&\hskip4cm=(x_i+\ov{x}_i-x_j-\ov{x}_j) h^{\oo}_{m-1}(x_i,\ov{x}_i,\ldots,x_{j},\ov{x}_{j},1|\a)\,; \label{eqn-oo-hra-xixj} \\
	&h^{\eo}_m(x_i,\ov{x}_i,\ldots,x_{j-1},\ov{x}_{j-1}|\a)-h^{\eo}_m(x_{i+1},\ov{x}_{i+1},\ldots,x_{j},\ov{x}_{j}|\a)\cr
	&\hskip4cm=(x_i+\ov{x}_i-x_j-\ov{x}_j) h^{\eo}_{m-1}(x_i,\ov{x}_i,\ldots,x_{j},\ov{x}_{j}|\a)\,; \label{eqn-eo-hra-xixj} 
\end{align}
\end{Lemma}

\noindent{\bf Proof}: First it should be noted that all these identities are trivially true for $m<0$ and for $m=0$ since
each $h_{m}$ reduces to either $0$ or $1$, with each $h_{m-1}$ reducing to $0$.
For $m>0$, in the simplest case 
\begin{align}
     &h^{\gl}_m(x_i,\ldots x_{j-1}|\a)-h^{\gl}_m(x_{i+1},\ldots,x_j|\a)\cr
		 &=[t^m]\ \left((1-tx_j)-(1-tx_i)\right) \ \prod_{\ell=i}^j \frac{1}{1-tx_\ell} \ \prod_{k=1}^{m+j-i-1}(1+ta_k)\cr
	   &= (x_i-x_j)\ [t^{m-1}]\  \prod_{\ell=i}^j \frac{1}{1-tx_\ell} \ \prod_{k=1}^{(m-1)+j-i}(1+ta_k)\cr
		 &= (x_i-x_j)\ h^{\gl}_{m-1}(x_i,\ldots,x_j|\a)\,. \label{eqn-gl-hrHa-xixj}
\end{align}

The other three cases are essentially the same, and are illustrated by the symplectic case:
\begin{align}
  &h^{\sp}_m(x_i,\ov{x}_i,\ldots,x_{j-1},\ov{x}_{j-1}|\a)-h^{\sp}_m(x_{i+1},\ov{x}_{i+1},\ldots,x_{j},\ov{x}_{j},|\a)\cr
&= [t^m]\ \left( (1\!-\!tx_j)(1\!-\!t\ov{x}_j)\!-\!(1\!-\!tx_i)(1\!-\!t\ov{x}_i)\right)\ \prod_{\ell=i}^j \frac{1}{(1\!-\!tx_\ell)(1\!-\!t\ov{x}_\ell)} 
          \!\!\prod_{k=1}^{m+j-i-1}\!\!(1+ta_k)\cr
		 &= (x_i\!+\!\ov{x}_i\!-\!x_j\!-\!\ov{x}_j)\ [t^{m-1}]\ \prod_{\ell=i}^j \frac{1}{(1\!-\!tx_\ell)(1\!-\!t\ov{x}_\ell)} 
		       \!\!\prod_{k=1}^{(m-1)+j-i}\!\!(1+ta_k)\cr 
	   &= (x_i+\ov{x}_i-x_j-\ov{x}_j)\ h^{\sp}_{m-1}(x_i,\ov{x}_i,\ldots,x_{j},\ov{x}_{j},|\a)\,. \label{eqn-sp-hrHa-xixj}
\end{align}

Exactly the same procedure applies to all the remaining cases, except for the even dimensional orthogonal group case with $m>0$ and $j=i+1$,
which is covered by the following argument:
\begin{align} 
   &h^{\eo}_{m}(x_i,\ov{x}_i|\a)-h^{\eo}_{m}(x_{i+1},\ov{x}_{i+1}|\a)\cr
   &= [t^m]\ \left(\frac{1}{1-tx_i}+\frac{1}{1-t\ov{x}_i}-\frac{1}{1-tx_{i+1}}-\frac{1}{1-t\ov{x}_{i+1}}\right) \prod_{k=1}^{m} (1+ta_k)\cr
   &= [t^{m-1}]\ \left( \frac{(1-t^2)(x_i+\ov{x}_i-x_{i+1}-\ov{x}_{i+1})}{(1-tx_i)(1-t\ov{x}_i)(1-tx_{i+1})(1-t\ov{x}_{i+1})}\right) \prod_{k=1}^m (1+ta_k)\cr
   &= (x_i+\ov{x}_i-x_{i+1}-\ov{x}_{i+1}) h^{\eo}_{m-1}(x_i,\ov{x}_i,x_{i+1},\ov{x}_{i+1}|\a)\,.  \label{eqn-eo-hrHa-xixj}
\end{align}
\qed

Now we are a position to state and prove the following result:

\begin{Theorem} [Flagged Jacobi-Trudi identity]\label{thm-fJT}
Given $\x=(x_1,x_2,\ldots,x_n)$ and $\ovx=(\ov{x}_1,\ov{x}_2,\ldots,\ov{x}_n)$, 
let $\x^{(i)}=(x_i,x_{i+1},\dots,x_n)$ and $\ovx^{(i)}=(\ov{x}_i,\ov{x}_{i+1},\ldots,\ov{x}_n)$. Then for any
partition $\lambda=(\lambda_1,\lambda_2,\ldots,\lambda_n)$ and any $\a=(a_1,a_2,\ldots)$ 
we have 
\begin{align}
\gl_\lambda(\x|\a) 
&= \left|\, h^{\gl}_{\lambda_j-j+i}(\x^{(i)}|\a) \,\right|\,; \label{eqn-gl-fJT} \\						
sp_\lambda(\x,\ovx|\a) 
&= \left|\, h^{\sp}_{\lambda_j-j+i}(\x^{(i)},\ovx^{(i)}|\a) \,\right| \,; \label{eqn-sp-fJT}\\
so_\lambda(\x,\ovx,1|\a)
&=  \left|\, h^{\oo}_{\lambda_j-j+i}(\x^{(i)},\ovx^{(i)},1|\a) \,\right|\,; \label{eqn-oo-fJT}\\
o_\lambda(\x,\ovx|\a)
&=  \left|\, h^{\eo}_{\lambda_j-j+i}(\x^{(i)},\ovx^{(i)}|\a) \,\right| \,, \label{eqn-eo-fJT}																	
\end{align}
\end{Theorem}

\noindent{\bf Proof}:
In the factorial Schur function $gl(n)$ case, we may proceed by manipulating the determinant in the numerator of the factorial Schur function formula 
taken from (\ref{eqn-gl-hm}). Subtracting row $(i+1)$ from row $i$ for $i=1,2,\ldots,n-1$ and applying (\ref{eqn-gl-hra-xixj}), 
then repeating the process for $i=1,2,\ldots,n-2$ and so on, yields
\begin{align}\label{eqn-gl-hlambda}
 \left|\, h^{\gl}_{\lambda_j+n-j}(x_i|\a) \,\right| 
&= \left|\,\begin{array}{c} h^{\gl}_{\lambda_j+n-j}(x_i|\a) - h^{\gl}_{\lambda_j+n-j}(x_{i+1}|\a) \cr
                                                                                  h^{\gl}_{\lambda_j+n-j}(x_n|\a) \cr 
																																									 \end{array}    \,\right| \cr
& = \prod_{i=1}^{n-1} (x_i-x_{i+1}) \ \left|\,\begin{array}{c} h^{\gl}_{\lambda_j+n-j-1}(x_i,x_{i+1}|\a) \cr
                                                                            h^{\gl}_{\lambda_j+n-j}(x_n|\a) \cr 
																																									 \end{array}  \,\right| \cr  
&= \prod_{i=1}^{n-1} (x_i-x_{i+1}) \ \left|\,\begin{array}{c} h^{\gl}_{\lambda_j+n-j-1}(x_i,x_{i+1}|\a)-h^{\gl}_{\lambda_j+n-j-1}(x_{i+1},x_{i+2} |\a) \cr
                                                                          h^{\gl}_{\lambda_j+n-j-1}(x_{n-1},x_n|\a) \cr 
																																						h^{\gl}_{\lambda_j+n-j}(x_n|\a) \cr 
    																																							 \end{array}  \,\right| \cr																																
& = \prod_{i=1}^{n-1} (x_i-x_{i+1})\prod_{i=1}^{n-2} (x_i-x_{i+2}) \ \left|\,\begin{array}{c} h^{\gl}_{\lambda_j+n-j-2}(x_i,x_{i+1},x_{i+2}|\a) \cr
                                                                           h^{\gl}_{\lambda_j+n-j-1}(x_{n-1},x_n|\a) \cr 
                                                                            h^{\gl}_{\lambda_j+n-j}(x_n|\a) \cr 
																																									 \end{array}  \,\right| \cr  
& = ~~\cdots~~ =   \prod_{1\leq i<j\leq n} (x_i-x_j)\ \left|\, h^{\gl}_{\lambda_j+n-j-(n-i)}(x_i,x_{i+1},\ldots,x_n|\a) \,\right| \cr
& = \prod_{1\leq i<j\leq n} (x_i-x_j)\ \left|\, h^{\gl}_{\lambda_j-j+i}(\x^{(i)}|\a) \,\right|\,. 
\end{align}
In the special case $\lambda=(0)$ we recover the denominator identity
\begin{equation}\label{eqn-gl-h0}
         \left|\, h^{\gl}_{n-j}(x_i|\a) \,\right|=\prod_{1\leq i<j\leq n} (x_i-x_j)\ \left|\, h^{\gl}_{-j+i}(\x^{(i)}|\a) \,\right|=\prod_{1\leq i<j\leq n} (x_i-x_j)
\end{equation} 
since the determinant $\left|\, h^{\gl}_{-j+i}(\x^{(i)}|\a) \,\right|$ is lower-triangular with all its diagonal elements equal to $h^{\gl}_0(\x|\a)=1$.
Taking the ratio of these two formulae implies that $gl_\lambda(\x|\a) =  \left|\, h^{\gl}_{\lambda_j-j+i}(\x^{(i)}|\a) \,\right|$, 
as required in (\ref{eqn-gl-fJT}).

The other three cases are almost identical. Taking the symplectic case as an example, and proceeding in exactly the same way as in the 
factorial Schur function case, one arrives at
\begin{align}
\left|\, h^{\sp}_{\lambda_j+n-j}(x_i|\a) \,\right| 
& = \prod_{i=1}^{n-1} (x_i+\ov{x}_i-x_{i+1}-\ov{x}_{i+1}) \ 
          \left|\,\begin{array}{c} h^{\sp}_{\lambda_j+n-j-1}(x_i,\ov{x}_i,x_{i+1},\ov{x}_{i+1}|\a) \cr
                                   h^{\sp}_{\lambda_j+n-j}(x_n,\ov{x}_n|\a) \cr 
									\end{array}  \,\right| \cr  
& = \prod_{i=1}^{n-1} (x_i+\ov{x}_i-x_{i+1}-\ov{x}_{i+1})\prod_{i=1}^{n-2} (x_i+\ov{x}_i-x_{i+2}-\ov{x}_{i+2}) \cr
&~~~~~~~~~~~~\ \left|\,\begin{array}{c} 
                      h^{\sp}_{\lambda_j+n-j-2}(x_i,\ov{x}_i,x_{i+1},\ov{x}_{i+1},x_{i+2},\ov{x}_{i+2}|\a) \cr
                      h^{\sp}_{\lambda_j+n-j-1}(x_{n-1},\ov{x}_{n-1},x_n,\ov{x}_n|\a) \cr 
                      h^{\sp}_{\lambda_j+n-j}(x_n,\ov{x}_n|\a) \cr 
										\end{array}  \,\right| \cr  
& = ~~\cdots~~  \cr   
& = \prod_{1\leq i<j\leq n} (x_i+\ov{x}_i-x_{j}-\ov{x}_{j})\ \left|\, h^{\sp}_{\lambda_j-j+i}(\x^{(i)},\ovx^{(i)}|\a) \,\right|\,. \label{eqn-sp-hlambda}
\end{align}
Applying this with $\lambda=(0)$ and taking the ratio of the these results gives (\ref{eqn-sp-fJT}), as required. The other two results follow in
exactly the same way.
\qed

\section{Explicit formulae in the case of one part partitions}
\label{sec-hma}

As a consequence of Theorem~\ref{thm-fJT} it should be noted that we have
\begin{Corollary}
In the special case $\lambda=(m,0,\ldots,0)$ we have
\begin{align}
\gl_{(m)}(\x|\a) &= h^{\gl}_{m}(\x|\a) \,;  &	so_{(m)}(\x,\ovx,1|\a) &= h^{\oo}_{m}(\x,\ovx,1|\a) \,; \cr
sp_{(m)}(\x,\ovx|\a) &= h^{\sp}_{m}(\x,\ovx|\a)  \,; & o_{(m)}(\x,\ovx|\a) &= h^{\eo}_{m}(\x,\ovx|\a)  \,. \nonumber
\end{align}
\end{Corollary}

\noindent{\bf Proof}: 
On setting $\lambda=(m,0,\ldots,0)$ the flagged Jacobi-Trudi determinants in (\ref{eqn-gl-fJT})-(\ref{eqn-eo-fJT}) 
are reduced to lower-triangular form since each $h_{-j+i}=0$ for $i<j$.  Moreover for $i>1$ the diagonal entries are
all of the form $h_0=1$, while the $(1,1)$ entry is just the appropriate $h_{m}$ with $\x^{(1)}=\x$ and $\ovx^{(1)}=\ovx$. 
\qed

Factorial characters in the one-part partition case $\lambda=(m,0,\ldots,0)$ with $m>0$ 
may then be evaluated directly from the generating function 
formulae of Definition~\ref{def-hma}. 
In the Schur function, or general linear case we have
\begin{Lemma}\label{lem-hm}
For $\x=(x_1,x_2,\ldots,x_n)$, $\a=(a_1,a_2,\ldots)$ and $m>0$
\begin{equation}\label{eqn-gl-hmxa}
  h^{\gl}_m(\x|\a)= \sum_{1\leq i_1\leq i_2\leq \cdots\leq i_m\leq n} (x_{i_1}+a_{i_1})(x_{i_2}+a_{i_2+1})\cdots(x_{i_m}+a_{i_m+m-1})\,.
\end{equation}
\end{Lemma}

\noindent{\bf Proof}:
First let the right hand side of (\ref{eqn-gl-hmxa}) be represented by $f_m(\x|\a)$, that is to say
\begin{equation}\label{eqn-gl-fmxa}
  f_m(\x|\a)= \sum_{1\leq i_1\leq i_2\leq \cdots\leq i_m\leq n} (x_{i_1}+a_{i_1})(x_{i_2}+a_{i_2+1})\cdots(x_{i_m}+a_{i_m+m-1})\,.
\end{equation}
It follows that for $m=1$
\begin{equation}\label{eqn-f1}
  f_1(\x|\a)=\sum_{i=1}^n (x_i+a_i)\,, \nonumber
\end{equation}
and for $m>1$
\begin{equation}\label{eqn-fm}
  f_m(\x|\a)= (x_1+a_1)f_{m-1}(\x|\tau\a)+f_m(\x^{(2)}|\tau\a)\,, \nonumber
\end{equation}
as can be seen by separating out the terms involving $x_1$. 

It remains only to show that $h^{\gl}_m(\x|\a)$ as defined by (\ref{eqn-gl-hma})
satisfies the same $m=1$ condition and the same recurrence relation. 
Expanding in powers of $t$ in the case $m=1$ immediately gives
\begin{equation}
    h^{\gl}_1(\x|\a) = [t]\ \prod_{i=1}^n \frac{1}{1-tx_i}\  \prod_{k=1}^{n}(1+ta_k) = \sum_{i=1}^n (x_i+a_i)\,, \nonumber
\end{equation}
as required.
While for $m>1$ we have
\begin{align}
  h^{\gl}_m(\x|\a) &= [t^m]\ \prod_{i=1}^n \frac{1}{1-tx_i}\  \prod_{k=1}^{m+n-1}(1+ta_k)\cr
	               &= [t^m]\ \left( \frac{1+ta_1}{1-tx_1}\right)\ \prod_{i=2}^{n} \frac{1}{1-tx_i}\  \prod_{k=2}^{m+n-1}(1+ta_k)\cr
								 &= [t^m]\ \left( 1+ \frac{t(x_1+a_1)}{1-tx_1}\right)\ \prod_{i=2}^{n} \frac{1}{1-tx_i}\  \prod_{k=2}^{m+n-1}(1+ta_k)\cr
								 &= (x_1+a_1) h^{\gl}_{m-1}(\x|\tau\a)+h^{\gl}_m(\x^{(2)}|\tau\a)\,, \label{eqn-gl-hmx1}
\end{align}
again precisely as required.
\qed

This result can be exploited in the symplectic case, where it might be noted that if we introduce dummy
parameters $a_\ell=0$ for $\ell\leq0$ we have
\begin{align}
  &h^{\sp}_m(\x,\ovx|\a) = [t^m]\ \prod_{i=1}^n \frac{1}{(1-tx_i)(1-t\ov{x}_i)}\  \prod_{k=1}^{m+n-1}(1+ta_k)\cr
	                         &= [t^m]\ \prod_{i=1}^n \frac{1}{(1-tx_i)(1-t\ov{x}_i)}\  \prod_{k=1-n}^{(m+2n-1)-n}(1+ta_k)
													 = h^{\gl}_m(\x,\ovx\,|\,\tau^{-n}\a)\,. \label{eqn-sphm-tau-n}
\end{align}
It follows that
\begin{equation}\label{eqn-sp-hm-z}
  h^{\sp}_m(\x,\ovx|\a)= \sum_{1\leq i_1\leq i_2\leq \cdots\leq i_m\leq 2n} (z_{i_1}+a_{i_1-n})(z_{i_2}+a_{i_2-n+1})\cdots(z_{i_m}+a_{i_m-n+m-1})\,.
\end{equation}
where 
\begin{equation}\label{eqn-sp-zwgts}
       z_{i_j}+a_{i_j-n+j-1}=\begin{cases} x_k+a_{2k-n+j-2}&\hbox{if $i_j=2k-1$;}\cr
			                                     \ov{x}_k+a_{2k-n+j-1}&\hbox{if $i_j=2k$,}\cr
																					\end{cases} ~~~~\hbox{with $a_\ell=0$ if $\ell\leq 0$}.
\end{equation}

Turning next to the odd orthogonal case, if we again introduce $a_\ell=0$ for $\ell\leq0$, we have
\begin{align}
  h^{\oo}_m(\x,\ovx,1|\a) &= [t^m]\ (1+t)\ \prod_{i=1}^n \frac{1}{(1-tx_i)(1-t\ov{x}_i)}\  \prod_{k=1}^{m+n-1}(1+ta_k)\cr
	                         &= [t^m]\ \left( \frac{1+t}{1+ta_{m+n}}\right) \prod_{i=1}^n \frac{1}{(1-tx_i)(1-t\ov{x}_i)}\  \prod_{k=1}^{m+n}(1+ta_k)\cr
													&= [t^m]\ \left(1 + \frac{t(1-a_{m+n})}{1+ta_{m+n}}\right) \prod_{i=1}^n \frac{1}{(1-tx_i)(1-t\ov{x}_i)}\  \prod_{k=1}^{m+n}(1+ta_k)\cr
													&= [t^m]\ \prod_{i=1}^n \frac{1}{(1-tx_i)(1-t\ov{x}_i)}\  \prod_{k=1+(1-n)}^{(m+2n-1)+(1-n)}(1+ta_k)\cr
													 &\ + (1-a_{m+n})\ [t^{m-1}]\ \prod_{i=1}^n \frac{1}{(1-tx_i)(1-t\ov{x}_i)}\  \prod_{k=1+(1-n)}^{(m-1+2n-1)+(1-n)}(1+ta_k)\cr
													&= h^{\gl}_m(\x,\ovx\,|\,\tau^{1-n} \a) + (1-a_{m+n})\ h^{\gl}_{m-1}(\x,\ovx\,|\,\tau^{1-n}\a)\,. \label{eqn-oohm-tau1-n}
\end{align}
It follows that
\begin{align}
  &h^{\oo}_m(\x,\ovx,1|\a)= \!\!\!\sum_{1\leq i_1\leq i_2\leq \cdots\leq i_m\leq 2n}\!\!\! (z_{i_1}+a_{i_1-n+1})(z_{i_2}+a_{i_2-n+2})\cdots(z_{i_m}+a_{i_m-n+m})\cr
  &~~~~~+\!\!\!\!\sum_{1\leq i_1\leq i_2\leq \cdots\leq i_{m-1}\leq 2n}\!\!\!\!(z_{i_1}+a_{i_1-n+1})\cdots(z_{i_{m-1}}+a_{i_{m-1}-n+m-1})(1-a_{m+n})\,, \label{eqn-oohm-xovx1}
\end{align}
where 
\begin{equation}\label{eqn-oo-zwgts}
       z_{i_j}+a_{i_j-n+j}=\begin{cases} x_k+a_{2k-n+j-1}&\hbox{if $i_j=2k-1$;}\cr
			                                     \ov{x}_k+a_{2k-n+j}&\hbox{if $i_j=2k$,}\cr
																					\end{cases} ~~~~\hbox{with $a_\ell=0$ if $\ell\leq 0$}.
\end{equation}
It should be noted that this expression must be independent of $a_{m+n}$, since this was introduced as a dummy parameter
in the form $(1+ta_{m+n})/(1+ta_{m+n})$. This independence can be verified by noting the cancellation that takes place 
between terms involving $(\ov{x}_n+a_{m+n})$ and $(1-a_{m+n})$ in (\ref{eqn-oohm-xovx1}). 

The case of the even orthogonal group is more difficult in part because of the distinction between the cases $n=1$ and $n>1$
that is already evident in (\ref{eqn-eo-hma}). In fact the case $n=1$ and $m>0$ is rather easy:
\begin{align}
h^{\eo}_m(x_1,\ov{x}_1|\a) &=  [t^m]\ \left(\frac{1}{1-tx_1}+\frac{1}{1-t\ov{x}_1} \right)\  \prod_{k=1}^{m}(1+ta_k)\cr
                               &=  h^{\gl}_m(x_1|\a) + h^{\gl}_m(\ov{x}_1|\a)
														    = (x_1|\a)^m+(\ov{x}_1|\a)^m\,. \label{eqn-eohm-n1}
\end{align}

In the case $n>1$ and $m>0$, yet again with $a_\ell=0$ for $\ell\leq0$, we have
\begin{align}
&h^{\eo}_m(\x,\ovx|\a) = [t^m]\ (1-t^2)\ \prod_{i=1}^n \frac{1}{(1-tx_i)(1-t\ov{x}_i)}\ \prod_{j=1}^{m+n-1}(1+ta_j)\cr
      &= [t^m]\ \left(\frac{tx_1}{1-tx_1}+\frac{t\ov{x}_1}{1-t\ov{x}_1}+1\right)\ \prod_{i=2}^n \frac{1}{(1-tx_i)(1-t\ov{x}_i)}\  \prod_{j=1}^{m+n-1}(1+ta_j)\cr
	 & = x_1\,[t^{m-1}]\ \frac{1}{1-tx_1}\ \prod_{i=2}^n \frac{1}{(1-tx_i)(1-t\ov{x}_i)}\  \prod_{j=1+(2-n)}^{(m-1+2n-1-1)+(2-n)}(1+ta_j)\cr 
	 & + \ov{x}_1\,[t^{m-1}]\ \frac{1}{1-t\ov{x}_1}\ \prod_{i=2}^n \frac{1}{(1-tx_i)(1-t\ov{x}_i)}\  \prod_{j=1+(2-n)}^{(m-1+2n-1-1)+(2-n)}(1+ta_j)\cr 
	 & + [t^m]\ \prod_{i=2}^n \frac{1}{(1-tx_i)(1-t\ov{x}_i)}\  \prod_{j=1+(2-n)}^{(m+2n-2-1)+(2-n)}(1+ta_j)\cr
	 &=~~~~ (x_1+a_{2-n}) h^{\gl}_{m-1}(x_1,x_2,\ov{x}_2,\ldots,x_n,\ov{x}_n\,|\,\tau^{2-n}\a)\cr
	 &~~~~ + (\ov{x}_1+a_{2-n}) h^{\gl}_{m-1}(\ov{x}_1,x_2,\ov{x}_2,\ldots,x_n,\ov{x}_n\,|\,\tau^{2-n}\a)\cr
	 &~~~~ + h^{\gl}_m(x_2,\ov{x}_2,\ldots,x_n,\ov{x}_n\,|\,\tau^{2-n}\a)\,, \label{eqn-eohm-tau2-n}
\end{align} 
where we have used the fact that $a_{2-n}=0$ for $n>1$. 
It will be noted that this expression remains valid in the case $n=1$, since in this case it reduces
to (\ref{eqn-eohm-n1}) as can be seen by using the recurrence relation (\ref{eqn-gl-hmx1}) with $\x=(x_1)$ and with $\ovx=(\ov{x}_1)$.
Thus (\ref{eqn-eohm-tau2-n}) is valid for all $n\geq1$ and $m>0$.

It follows that
\begin{align}
  &h^{\eo}_m(\x,\ovx|\a)=\sum_{k=1}^m (x_1+a_{2-n})(x_1+a_{3-n})\cdots(x_1+a_{k+1-n})\cr
	&\times\!\!\!\sum_{3\leq i_{k+1}\leq i_{k+2}\leq \cdots\leq i_{m}\leq 2n}\!\!\! (z_{i_{k+1}}+a_{i_{k+1}-n+k})(z_{i_{k+2}}+a_{i_{k+2}-n+k+1})\cdots(z_{i_m}+a_{i_m-n+m-1})\cr
  &~~~~~+\sum_{k=1}^m (\ov{x}_1+a_{2-n})(\ov{x}_1+a_{3-n})\cdots(\ov{x}_1+a_{k+1-n})\cr
	&\times\!\!\!\sum_{3\leq i_{k+1}\leq i_{k+2}\leq \cdots\leq i_{m}\leq 2n}\!\!\!\!(z_{i_{k+1}}+a_{i_{k+1}-n+k})(z_{i_{k+2}}+a_{i_{k+2}-n+k+1})\cdots(z_{i_m}+a_{i_m-n+m-1}) \cr
	&~~~~~+\!\!\!\sum_{3\leq i_1\leq i_2\leq \cdots\leq i_m\leq 2n}\!\!\! (z_{i_1}+a_{i_1-n})(z_{i_2}+a_{i_2-n+1})\cdots(z_{i_m}+a_{i_m-n+m-1})\,,\cr \label{eqn-eohm-xovx1}
\end{align}
where 
\begin{equation}\label{eqn-eo-zwgts}
       z_{i_j}+a_{i_j-n+j}=\begin{cases} x_k+a_{2k-n+j-2+\delta_{i_j,1}}&\hbox{if $i_j=2k-1$;}\cr
			                                     \ov{x}_k+a_{2k-n+j-1}&\hbox{if $i_j=2k$,}\cr
																					\end{cases} ~~~~\hbox{with $a_\ell=0$ if $\ell\leq 0$}.
\end{equation}

\medskip


\section{Combinatorial realisation of factorial characters}
\label{sec-tableaux}

The significance of these results is that they offer various lattice path models of each of the relevant
one-part partition factorial characters. By making use of $n$-tuples of such lattice paths in the interpretation
of the flagged Jacobi-Trudi identities of Theorem~\ref{thm-fJT} one arrives at a non-intersecting
lattice path model of factorial characters specified by any partition $\lambda$ of length $\ell(\lambda)\leq n$. 
This leads inexorably to a further realisation of factorial characters in terms of certain
appropriately weighted tableaux. The tableaux themselves include those already associated with
Schur functions, symplectic group characters, and both odd and even orthogonal group characters in the classical non-factorial case.

Restricting our attention to fixed $n$ and partitions
$\lambda=(\lambda_1,\lambda_2,\ldots,\lambda_n)$ of length $\ell(\lambda)\leq n$,
each such partition defines a Young diagram $F^\lambda$
consisting of $|\lambda|=\lambda_1+\lambda_2+\cdots+\lambda_n$ boxes arranged in $\ell(\lambda)$ rows of lengths $\lambda_i$.
for $i=1,2,\ldots,\ell(\lambda)$. We adopt the English convention as used by Macdonald~\cite{Mac95} whereby the rows are
left-adjusted to a vertical line and are weakly decreasing in length from top to bottom.
For example in the case $n=6$ and $\lambda=(4,4,3,1,0,0)$, we have $\ell(\lambda)=4$, $|\lambda|=12$ and
$$
F^{4431}\ =\ 
\YT{0.2in}{0.2in}{}{
 {{},{},{},{}},
 {{},{},{},{}},
 {{},{},{}},
 {{}},
}
$$
More precisely we define $F^\lambda=\{ (i,j)\,|\, 1\leq i\leq \ell(\lambda) ; 1\leq j\leq \lambda_i\}$ and refer
to $(i,j)$ as being the box in the $i$th row and $j$th column of $F^\lambda$. Assigning an entry $T_{ij}$
taken from some alphabet $I$ to each box $(i,j)$ of $F^\lambda$ in accordance with various rules gives rise to 
tableaux $T$ of shape $\lambda$ that may be used, as we shall see, to express both ordinary and factorial characters 
in a combinatorial manner. We deal here with the cases $GL(n,\C)$, $SO(2n+1,\C)$, $Sp(2n,\C)$ and $O(2n,\C)$.

\begin{Definition}\label{def-glTab} (Case $gl(n)$): 
Let ${\cal T}^{\gl}_\lambda$ be the set of all semistandard Young tableaux $T$ of shape $\lambda$ 
that are obtained by filling each box $(i,j)$ of $F^\lambda$ with an entry $T_{ij}$ from the alphabet
$$
     \{ 1<2<\cdots<n \}
$$
in such a way that:
{\bf (T1)}: entries weakly increase from left to right across rows;
{\bf (T2)}: entries weakly increase from top to bottom down columns;
{\bf (T3)}: no two identical entries appear in any column.
\end{Definition}

\begin{Definition}\label{def-spTab} (Case $sp(2n)$): 
Let ${\cal T}^{\sp}_\lambda$ be the set of all symplectic tableau $T$ of shape $\lambda$
that are obtained by filling each box $(i,j)$ of $F^\lambda$ with an entry $T_{ij}$ from
the alphabet
$$
    \{1<\ov{1}<2<\ov{2}<\cdots<n<\ov{n}\}
$$
in such a way that the entries satisfy {\bf (T1)-(T3)}, along with
{\bf (T4)}: neither $k$ nor $\ov{k}$ appear lower than the $k$th row.

\end{Definition}
\begin{Definition}\label{def-ooTab} (Case $so(2n+1)$): 
Let ${\cal T}^{\oo}_\lambda$ be the set of all odd orthogonal tableaux $T$ of shape $\lambda$
that are obtained by filling each box $(i,j)$ of $F^\lambda$ with an entry $T_{ij}$ from
the alphabet
$$
    \{1<\ov{1}<2<\ov{2}<\cdots<n<\ov{n}<0\}
$$
in such a way that {\bf (T1)} and {\bf (T2)} are satisfied by all entries, along with 
{\bf (T3)} and {\bf (T4)} applying to non-zero entries, and  
{\bf (T5)}: at most one entry $0$ appears in any row.
\end{Definition}

\begin{Definition}\label{def-eoTab} (Case $o(2n)$): 
Let ${\cal T}^{\eo}_\lambda$ be the set of all even orthogonal tableau $T$ of shape $\lambda$
that are obtained by filling each box $(i,j)$ of $F^\lambda$ with an entry $T_{ij}$ from
the alphabet
$$
 \{1<\ov{1}<2<\ov{2}<\cdots<n<\ov{n}\}
$$
in such a way that the entries satisfy {\bf (T1)-(T4)} 
together with: 
{\bf (T6)} if an entry $k$ appears in the $k$th row then any entry $\ov{k}$ in the same row must
be covered by an entry $k$ appearing immediately above it in the preceding row.

\end{Definition}

These definitions are exemplified for $gl(4)$,  $sp(8)$, $so(9)$ and $o(8)$ 
as below from left to right:
\begin{equation}
\YT{0.2in}{0.2in}{}{
 {{1},{1},{2},{4}},
 {{2},{3},{3}},
 {{4},{4},{4}},
} 
\quad
\YT{0.2in}{0.2in}{}{
 {{1},{\ov1},{2},{\ov4}},
 {{\ov3},{4},{4}},
 {{4},{\ov4},{\ov4}},
}
\quad
\YT{0.2in}{0.2in}{}{
 {{1},{\ov1},{2},{\ov4}},
 {{3},{4},{0}},
 {{4},{\ov4},{0}},
}
\quad 
\YT{0.2in}{0.2in}{}{
 {{2},{2},{\ov2},{\ov2},{4}},
 {{\ov2},{3},{3},{\ov3},{\ov4}},
 {{3},{4},{4},{\ov4}},
 {{4},{\ov4},{\ov4}},
}
\end{equation}

These definitions allow us to provide combinatorial expressions for factorial characters as follows:
\begin{Theorem}\label{thm-Twgts}
For each $\g$ and $\z$ as tabulated below, 
\begin{equation}\label{eqn-factchar-tab}
   \g_\lambda(\z|\a) = \sum_{T\in{\cal T}^{\g}_\lambda}\ 2^{\zeta(T)}\prod_{(i,j)\in F^\lambda} \wgt(T_{ij})\,.
\end{equation}
with
\begin{equation}\label{eqn-Twgts}
\begin{array}{|l|l|l|l|l|l|}
\hline
\g_\lambda(\z|\a)&{\cal T}^\g&\zeta(T)&T_{ij}&\wgt(T_{ij})&\cr
\hline
gl_\lambda(\x|\a)&T\in{\cal T}&0&k &x_k+a_{k+j-i}&\cr
\hline
so_\lambda(\x,\ovx,1|\a)&T\in{\cal T}^{\oo}&0&k&x_k+a_{2k-n+j-i}&a_m=0\hbox{~~for~~}m\leq0\cr
                    &&&\ov{k}&\ov{x}_k+a_{2k+1-n+j-i}&\cr
										&&&0&1-a_{n+1+j-i}&\cr
\hline
sp_\lambda(\x,\ovx|\a)&T\in{\cal T}^{\sp}&0&k&x_k+a_{2k-1-n+j-i}&a_m=0\hbox{~~for~~}m\leq0\cr
                    &&&\ov{k}&\ov{x}_k+a_{2k-n+j-i}&\cr
\hline
o_\lambda(\x,\ovx|\a)&T\in{\cal T}^{\eo}&\sum_k\zeta_{k\ov{k}}&k&x_k+a_{2k-1-n+j-i+\delta_{ik}}&a_m=0\hbox{~~for~~}m\leq0\cr
                   &&&\ov{k}&\ov{x}_k+a_{2k-n+j-i}&\cr										
\hline					
\end{array}
\end{equation}
where $\zeta_{k\ov{k}}$ is $1$ if $T_{k-1,1}=k$ and $T_{k,1}=\ov{k}$, and is $0$ otherwise. 
\end{Theorem}

\bigskip
\noindent{\bf Proof}: {\bf Case $gl(n)$}.
The first of the above results in (\ref{eqn-Twgts}), namely the Schur function $gl(n)$ case, has been derived elsewhere~\cite{HK15b} 
by means of lattice path methods, but since it informs all the other cases, we repeat the argument here.
We adopt matrix coordinates $(k,\ell)$ for lattice points with $k=1,2,\ldots,n$ specifying row labels from
top to bottom, and $\ell=1,2,\ldots,\lambda_1+n$ specifying column labels from left to right. 
For example, if $n=4$ and $\lambda=(4,3,3,0)$ the lattice takes the form:
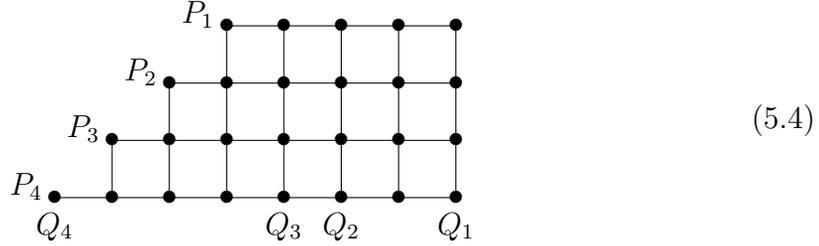
\begin{equation}\label{eqn-gln-grid}
\vcenter{\hbox{
\begin{tikzpicture}[x={(0in,-0.3in)},y={(0.3in,0in)}] 
\foreach \j in {4,...,8} \draw(1,\j)node{$\bullet$};
\foreach \j in {3,...,8} \draw(2,\j)node{$\bullet$};
\foreach \j in {2,...,8} \draw(3,\j)node{$\bullet$};
\foreach \j in {1,...,8} \draw(4,\j)node{$\bullet$};
\draw[-](1,4)to(1,8);
\draw[-](2,3)to(2,8);
\draw[-](3,2)to(3,8);
\draw[-](4,1)to(4,8);
\draw[-](2,3)to(4,3);
\draw[-](3,2)to(4,2);
\foreach \j in {4,...,8} \draw[-](1,\j)to(4,\j);
\draw(1-0.2,3.5)node{$P_1$};
\draw(2-0.2,2.5)node{$P_2$};
\draw(3-0.2,1.5)node{$P_3$};
\draw(4-0.2,0.5)node{$P_4$};
\draw(4.5,1)node{$Q_{4}$};  
\draw(4.5,5)node{$Q_{3}$};
\draw(4.5,6)node{$Q_{2}$};
\draw(4.5,8)node{$Q_{1}$};
\end{tikzpicture}
}}
\end{equation}
Each lattice path that we are interested in is a continuous path from some 
$P_i=(i,n-i+1)$ to some $Q_j=(n,n-j+1+\lambda_j)$ with $i,j\in\{1,2,\ldots,n\}$. 
Such a path consists of a sequence of horizontal or vertical edges
and is associated with a contribution to $h_{\lambda_j-j+i}(\x^{(i)}|\a)$ in the form of
a summand of (\ref{eqn-gl-hmxa}) with $m=\lambda_j-j+i$ and $\x$ replaced by $\x^{(i)}$. 
Taking into account the restriction of the alphabet from $\x$ to $\x^{(i)}$, the weight assigned to  
horizontal edge from $(k,\ell-1)$ to $(k,\ell)$ is $x_k+a_{k+\ell-n-1}$. 
Thanks to the Lindstr\"om-Gessel-Viennot theorem~\cite{Lin73,GV85,GV89} the only
surviving contributions to the determinantal expression for $s_\lambda(\x|\a)$ in the flagged factorial Jacobi-Trudi 
identity (\ref{eqn-gl-fJT}) are those corresponding to an $n$-tuple of non-intersecting lattice paths from $P_i$ to $Q_i$ for $i=1,2,\ldots,n$. 
Such $n$-tuples are easily seen~\cite{SW86} to be in bijective correspondence
with semistandard Young tableaux $T\in{\cal T}^{\gl}_\lambda$ of shape $\lambda$ as in Definition~\ref{def-glTab}, with the $j$th horizontal edge at level $k$
on the path from $P_i$ to $Q_i$ giving an entry $T_{ij}=k$ in $T$ for $i=1,2,\ldots,n$ and $j=1,2,\ldots,\lambda_i$. To complete the proof 
of Theorem~\ref{thm-Twgts} in the factorial $gl(n)$ case 
it only remains to note that the weight $\wgt(T_{ij})$ to be assigned to $T_{ij}$ is that of the edge from $(k,\ell-1)$ to $(k,\ell)$  
given by $x_k+a_{k+\ell-n-1}=x_k+a_{k+j-i}$ with $j=\ell-(n-i+1)$ since this is the number of horizontal steps from $P_i$ to column $\ell$
on the lattice path from $P_i$ to $Q_i$.
This is exemplified 
in the case $n=4$ and $\lambda=(4,3,3)$ by

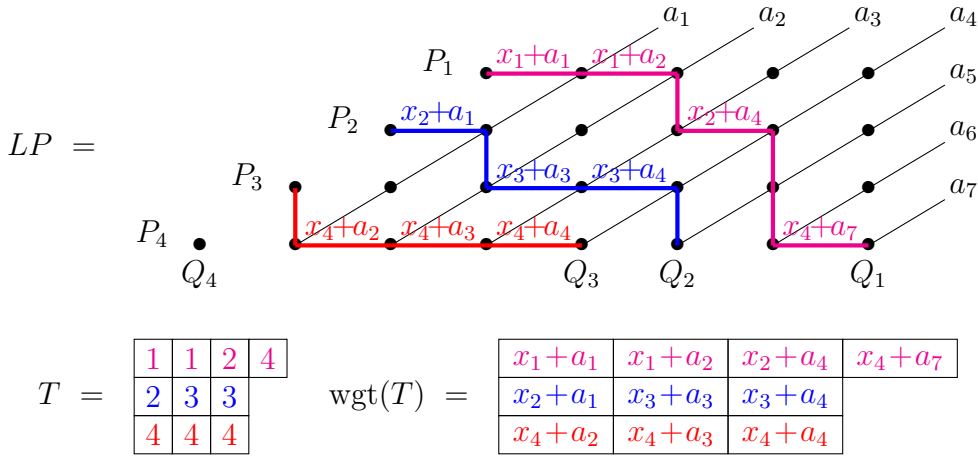
\begin{equation}\label{eqn-gl-LPT}
\begin{array}{c}
LP\ =\ 
\vcenter{\hbox{
\begin{tikzpicture}[x={(0in,-0.3in)},y={(0.5in,0in)}] 
\foreach \j in {4,...,8} \draw(1,\j)node{$\bullet$};
\foreach \j in {3,...,8} \draw(2,\j)node{$\bullet$};
\foreach \j in {2,...,8} \draw(3,\j)node{$\bullet$};
\foreach \j in {1,...,8} \draw(4,\j)node{$\bullet$};
\draw[-](0+0.2,6-0.2)to(4,2);
\draw[-](0+0.2,7-0.2)to(4,3);
\draw[-](0+0.2,8-0.2)to(4,4);
\draw[-](0+0.2,9-0.2)to(4,5);
\draw[-](1+0.2,9-0.2)to(4,6);
\draw[-](2+0.2,9-0.2)to(4,7);
\draw[-](3+0.2,9-0.2)to(4,8);
\draw(1-0.2,3.5)node{$P_1$};
\draw(2-0.2,2.5)node{$P_2$};
\draw(3-0.2,1.5)node{$P_3$};
\draw(4-0.2,0.5)node{$P_4$};
\draw(4.5,1)node{$Q_{4}$};  
\draw(4.5,5)node{$Q_{3}$};
\draw(4.5,6)node{$Q_{2}$};
\draw(4.5,8)node{$Q_{1}$};
\draw(0,6)node{$a_1$};
\draw(0,7)node{$a_2$};
\draw(0,8)node{$a_3$};
\draw(0,9)node{$a_4$};
\draw(1,9)node{$a_5$};
\draw(2,9)node{$a_6$};
\draw(3,9)node{$a_7$};
\draw[draw=magenta,ultra thick] (1,4)to(1,5); \draw(1-0.3,5-0.5)node{$\magenta{x_1\!\!+\!\!a_1}$};
\draw[draw=magenta,ultra thick] (1,5)to(1,6); \draw(1-0.3,6-0.5)node{$\magenta{x_1\!\!+\!\!a_2}$};
\draw[draw=magenta,ultra thick] (1,6)to(2,6); 
\draw[draw=magenta,ultra thick] (2,6)to(2,7); \draw(2-0.3,7-0.5)node{$\magenta{x_2\!\!+\!\!a_4}$};
\draw[draw=magenta,ultra thick] (2,7)to(4,7);
\draw[draw=magenta,ultra thick] (4,7)to(4,8); \draw(4-0.3,8-0.5)node{$\magenta{x_4\!\!+\!\!a_7}$};
\draw[draw=blue,ultra thick] (2,3)to(2,4); 
\draw[draw=blue,ultra thick] (2,4)to(3,4); \draw(2-0.3,4-0.5)node{$\blue{x_2\!\!+\!\!a_1}$};
\draw[draw=blue,ultra thick] (3,4)to(3,5); \draw(3-0.3,5-0.5)node{$\blue{x_3\!\!+\!\!a_3}$};
\draw[draw=blue,ultra thick] (3,5)to(3,6); \draw(3-0.3,6-0.5)node{$\blue{x_3\!\!+\!\!a_4}$};
\draw[draw=blue,ultra thick] (3,6)to(4,6);
\draw[draw=red,ultra thick] (3,2)to(4,2); 
\draw[draw=red,ultra thick] (4,2)to(4,3); \draw(4-0.3,3-0.5)node{$\red{x_4\!\!+\!\!a_2}$};
\draw[draw=red,ultra thick] (4,3)to(4,4); \draw(4-0.3,4-0.5)node{$\red{x_4\!\!+\!\!a_3}$}; 
\draw[draw=red,ultra thick] (4,4)to(4,5); \draw(4-0.3,5-0.5)node{$\red{x_4\!\!+\!\!a_4}$};
\end{tikzpicture}
}}\cr\cr
T\ = \ \YT{0.2in}{0.2in}{}{
 {\magenta{1},\magenta{1},\magenta{2},\magenta{4}},
 {\blue{2},\blue{3},\blue{3}},
 {\red{4},\red{4},\red{4}},
}
\quad
\wgt(T)\ =\ 
\wideYT{0.2in}{0.6in}{}{
 {\magenta{x_1\!+\!a_1},\magenta{x_1\!+\!a_2},\magenta{x_2\!+\!a_4},\magenta{x_4\!+\!a_7}},
 {\blue{x_2\!+\!a_1},\blue{x_3\!+\!a_3},\blue{x_3\!+\!a_4}},
 {\red{x_4\!+\!a_2},\red{x_4\!+\!a_3},\red{x_4\!+\!a_4}},
} \nonumber
\end{array}
\end{equation}

\bigskip
\noindent{\bf Case $sp(2n)$}. 
In the symplectic case, the use of (\ref{eqn-sphm-tau-n}) in (\ref{eqn-sp-fJT}) implies that
\begin{equation}
 sp_\lambda(\x,\ovx|\a)=\left|\, h^{\gl}_{\lambda_j-j+i}(\x^{(i)},\ovx^{(i)}\,|\,\tau^{i-1-n}\a) \,\right| 
\end{equation}
with $a_\ell=0$ for $\ell\leq0$,
the lattice path proof proceeds exactly as in the Schur function case
with $\a$ replaced by $\tau^{i-1-n}\a$, and with the alphabet extended to include both $x_k$ and $\ov{x}_k$ 
for $k=1,2,\ldots,n$ and the underlying lattice taking the typical form:
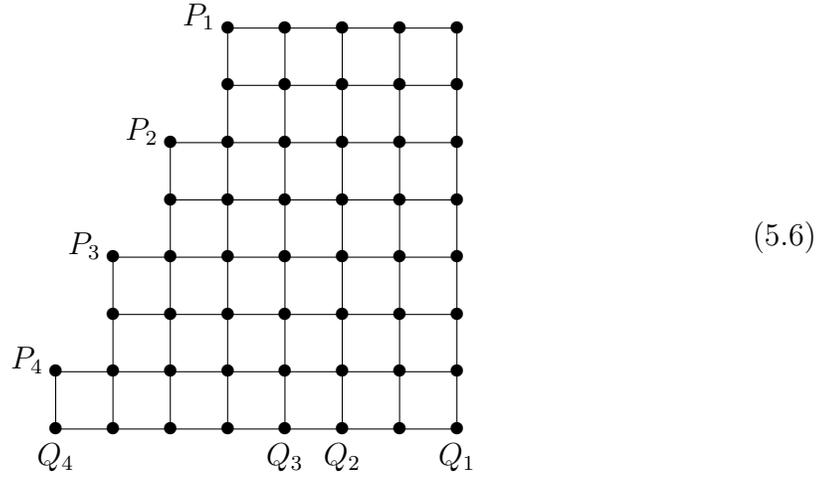
\begin{equation}\label{eqn-sp-grid}
\vcenter{\hbox{
\begin{tikzpicture}[x={(0in,-0.3in)},y={(0.3in,0in)}] 
\foreach \j in {4,...,8} \draw(1,\j)node{$\bullet$};
\foreach \j in {4,...,8} \draw(2,\j)node{$\bullet$};
\foreach \j in {3,...,8} \draw(3,\j)node{$\bullet$};
\foreach \j in {3,...,8} \draw(4,\j)node{$\bullet$};
\foreach \j in {2,...,8} \draw(5,\j)node{$\bullet$};
\foreach \j in {2,...,8} \draw(6,\j)node{$\bullet$};
\foreach \j in {1,...,8} \draw(7,\j)node{$\bullet$};
\foreach \j in {1,...,8} \draw(8,\j)node{$\bullet$};
\draw[-](1,4)to(1,8);
\draw[-](2,4)to(2,8);
\draw[-](3,3)to(3,8);
\draw[-](4,3)to(4,8);
\draw[-](5,2)to(5,8);
\draw[-](6,2)to(6,8);
\draw[-](7,1)to(7,8);
\draw[-](8,1)to(8,8);
\draw[-](7,1)to(8,1);
\draw[-](5,2)to(8,2);
\draw[-](3,3)to(8,3);
\foreach \j in {4,...,8} \draw[-](1,\j)to(8,\j);
\draw(1-0.2,3.5)node{$P_1$};
\draw(3-0.2,2.5)node{$P_2$};
\draw(5-0.2,1.5)node{$P_3$};
\draw(7-0.2,0.5)node{$P_4$};
\draw(8.5,1)node{$Q_{4}$};  
\draw(8.5,5)node{$Q_{3}$};
\draw(8.5,6)node{$Q_{2}$};
\draw(8.5,8)node{$Q_{1}$};
\end{tikzpicture}
}}
\end{equation}
The starting points are now $P_i=(2i-1,n-i+1)$ thereby ensuring that condition {\bf (T4)} is satisfied, 
and the end points are $Q_j=(2n,n-j+1+\lambda_j)$. Once again
it is only the $n$-tuples of non-intersecting lattice paths from $P_i$ to $Q_i$ that
contribute to $sp_\lambda(\x,\ovx|\a)$ and these are in bijective correspondence with
the symplectic tableaux of Definition~\ref{def-spTab} of shape $\lambda$ with entries from $\{1<\ov1<\cdots<n<\ov{n}\}$.
This is exemplified 
for $n=4$ and $\lambda=(4,3,3,0)$ by the following, where we have used the fact that $a_m=0$ for $m\leq0$ in the
tableau weights but not in the lattice path edge weights: 

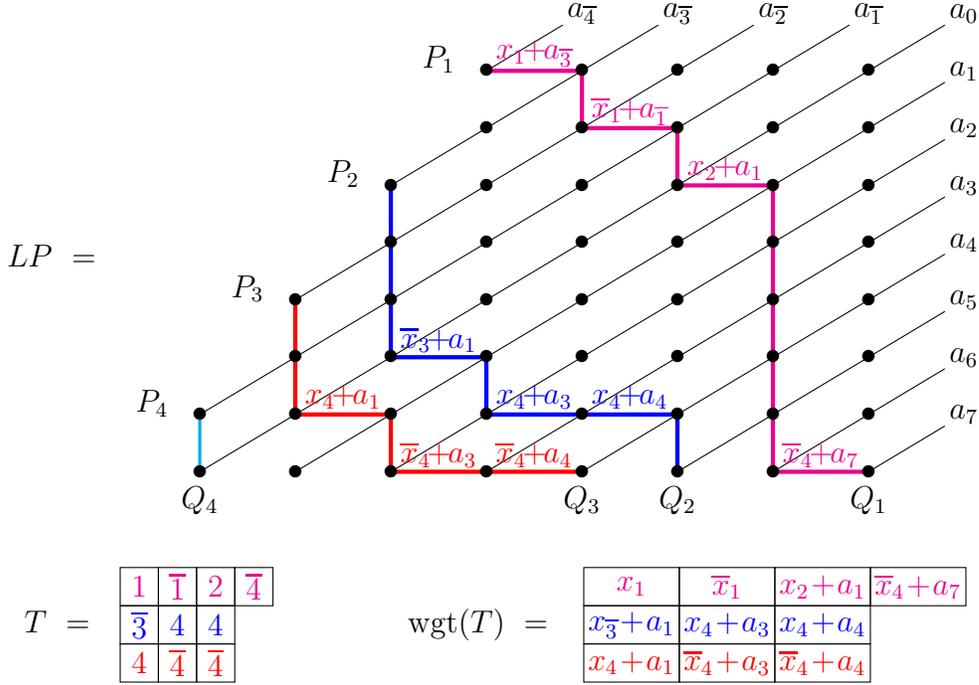
\begin{equation}\label{eqn-sp-LPTex}
\begin{array}{c}
LP\ = \ 
\vcenter{\hbox{
\begin{tikzpicture}[x={(0in,-0.3in)},y={(0.5in,0in)}] 
\draw[-](0+0.2,5-0.2)to(1,4);
\draw[-](0+0.2,6-0.2)to(3,3);
\draw[-](0+0.2,7-0.2)to(5,2);
\draw[-](0+0.2,8-0.2)to(7,1);
\foreach \i in {0,...,7} \draw[-,thin](\i+0.2,9-0.2)to(8,\i+1);
\draw(1-0.2,4-0.5)node{$P_1$};
\draw(3-0.2,3-0.5)node{$P_2$};
\draw(5-0.2,2-0.5)node{$P_3$};
\draw(7-0.2,1-0.5)node{$P_4$};
\draw(8+0.5,1)node{$Q_4$};
\draw(8+0.5,5)node{$Q_3$};
\draw(8+0.5,6)node{$Q_2$};
\draw(8+0.5,8)node{$Q_1$};
\draw(0,5)node{$a_{\ov4}$};
\draw(0,6)node{$a_{\ov3}$};
\draw(0,7)node{$a_{\ov2}$};
\draw(0,8)node{$a_{\ov1}$};
\foreach \i in {0,...,7} \draw(\i,9)node{$a_\i$};
\draw[draw=magenta,ultra thick] (1,4)to(1,5); \draw(1-0.3,5-0.5)node{$\magenta{x_1\!\!+\!\!a_{\ov3}}$};
\draw[draw=magenta,ultra thick] (1,5)to(2,5);
\draw[draw=magenta,ultra thick] (2,5)to(2,6); \draw(2-0.3,6-0.5)node{$\magenta{\ov{x}_1\!\!+\!\!a_{\ov1}}$};
\draw[draw=magenta,ultra thick] (2,6)to(3,6);
\draw[draw=magenta,ultra thick] (3,6)to(3,7); \draw(3-0.3,7-0.5)node{$\magenta{x_2\!\!+\!\!a_1}$};
\draw[draw=magenta,ultra thick] (3,7)to(8,7);
\draw[draw=magenta,ultra thick] (8,7)to(8,8); \draw(8-0.3,8-0.5)node{$\magenta{\ov{x}_4\!\!+\!\!a_7}$};
\draw[draw=blue,ultra thick] (3,3)to(6,3);
\draw[draw=blue,ultra thick] (6,3)to(6,4);  \draw(6-0.3,4-0.5)node{$\blue{\ov{x}_3\!\!+\!\!a_1}$};
\draw[draw=blue,ultra thick] (6,4)to(7,4);
\draw[draw=blue,ultra thick] (7,4)to(7,5);  \draw(7-0.3,5-0.5)node{$\blue{x_4\!\!+\!\!a_3}$};
\draw[draw=blue,ultra thick] (7,5)to(7,6);  \draw(7-0.3,6-0.5)node{$\blue{x_4\!\!+\!\!a_4}$};
\draw[draw=blue,ultra thick] (7,6)to(8,6);
\draw[draw=red,ultra thick] (5,2)to(7,2);
\draw[draw=red,ultra thick] (7,2)to(7,3);  \draw(7-0.3,3-0.5)node{$\red{x_4\!\!+\!\!a_1}$};
\draw[draw=red,ultra thick] (7,3)to(8,3);
\draw[draw=red,ultra thick] (8,3)to(8,4);  \draw(8-0.3,4-0.5)node{$\red{\ov{x}_4\!\!+\!\!a_3}$};
\draw[draw=red,ultra thick] (8,4)to(8,5);  \draw(8-0.3,5-0.5)node{$\red{\ov{x}_4\!\!+\!\!a_4}$};
\draw[draw=cyan,very thick] (7,1)to(8,1);
\foreach \j in {4,...,8} \draw(1,\j)node{$\bullet$};
\foreach \j in {4,...,8} \draw(2,\j)node{$\bullet$};
\foreach \j in {3,...,8} \draw(3,\j)node{$\bullet$};
\foreach \j in {3,...,8} \draw(4,\j)node{$\bullet$};
\foreach \j in {2,...,8} \draw(5,\j)node{$\bullet$};
\foreach \j in {2,...,8} \draw(6,\j)node{$\bullet$};
\foreach \j in {1,...,8} \draw(7,\j)node{$\bullet$};
\foreach \j in {1,...,8} \draw(8,\j)node{$\bullet$};
\end{tikzpicture}
}}\cr\cr
T\ = \
\YT{0.2in}{0.2in}{}{
 {\magenta{1},\magenta{\ov1},\magenta{2},\magenta{\ov4}},
 {\blue{\ov3},\blue{4},\blue{4}},
 {\red{4},\red{\ov4},\red{\ov4}},
}
\qquad\qquad
\wgt(T)\ = \ 
\YT{0.2in}{0.5in}{}{
 {\magenta{x_1},\magenta{\ov{x}_1},\magenta{{x_2}\!+\!{a_1}},\magenta{\ov{x}_4\!+\!{a_7}}},
 {\blue{x_{\ov3}\!+\!{a_1}},\blue{{x_4}\!+\!{a_3}},\blue{{x_4}\!+\!{a_4}}},
 {\red{{x_4}\!+\!{a_1}},\red{\ov{x}_4\!+\!{a_3}},\red{\ov{x}_4\!+\!{a_4}}},
} \nonumber
\end{array}
\end{equation}

\bigskip
\noindent{\bf Case $so(2n+1)$}.
In the odd orthogonal case, the use of (\ref{eqn-oohm-tau1-n}) in (\ref{eqn-oo-fJT}) which gives
\begin{equation}
so_\lambda(\x,\ovx,1|\a)=\left|\, h^{\gl}_{\lambda_j-j+i}(\x^{(i)},\ovx^{(i)}\,|\,\tau^{i-n} \a) 
           + (1-a_{\lambda_j-j+n+1})\ h^{\gl}_{\lambda_j-j+i-1}(\x^{(i)},\ovx^{(i)}\,|\,\tau^{i-n}\a) \,\right| \,.
\end{equation} 
As indicated in (\ref{eqn-oohm-xovx1}), this involves extending the alphabet to include not only both $x_k$ and $\ov{x}_k$ for $k=1,2,\ldots,n$, 
but also $1$, and $\a$ is replaced this time by $\tau^{i-n}\a$. The underlying grid typically takes the form:
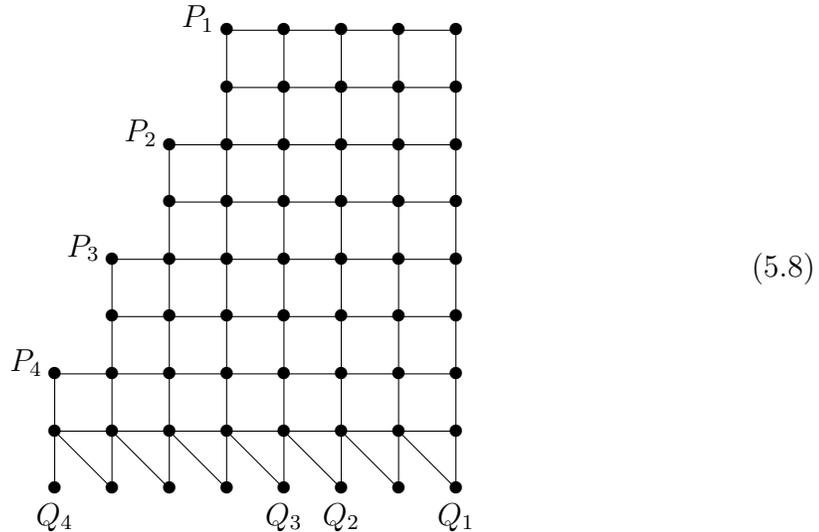
\begin{equation}\label{eqn-oo-grid}
\vcenter{\hbox{
\begin{tikzpicture}[x={(0in,-0.3in)},y={(0.3in,0in)}] 
\foreach \j in {4,...,8} \draw(1,\j)node{$\bullet$};
\foreach \j in {4,...,8} \draw(2,\j)node{$\bullet$};
\foreach \j in {3,...,8} \draw(3,\j)node{$\bullet$};
\foreach \j in {3,...,8} \draw(4,\j)node{$\bullet$};
\foreach \j in {2,...,8} \draw(5,\j)node{$\bullet$};
\foreach \j in {2,...,8} \draw(6,\j)node{$\bullet$};
\foreach \j in {1,...,8} \draw(7,\j)node{$\bullet$};
\foreach \j in {1,...,8} \draw(8,\j)node{$\bullet$};
\foreach \j in {1,...,8} \draw(9,\j)node{$\bullet$};
\draw[-](1,4)to(1,8);
\draw[-](2,4)to(2,8);
\draw[-](3,3)to(3,8);
\draw[-](4,3)to(4,8);
\draw[-](5,2)to(5,8);
\draw[-](6,2)to(6,8);
\draw[-](7,1)to(7,8);
\draw[-](8,1)to(8,8);
\draw[-](7,1)to(9,1);
\draw[-](5,2)to(9,2);
\draw[-](3,3)to(9,3);
\foreach \j in {4,...,8} \draw[-](1,\j)to(9,\j);
\foreach \j in {1,...,7} \draw[-](8,\j)to(9,\j+1);
\draw(1-0.2,3.5)node{$P_1$};
\draw(3-0.2,2.5)node{$P_2$};
\draw(5-0.2,1.5)node{$P_3$};
\draw(7-0.2,0.5)node{$P_4$};
\draw(9.5,1)node{$Q_{4}$};  
\draw(9.5,5)node{$Q_{3}$};
\draw(9.5,6)node{$Q_{2}$};
\draw(9.5,8)node{$Q_{1}$};
\end{tikzpicture}
}}
\end{equation}
The starting points are $P_i=(2i-1,n-i+1)$, 
ensuring, as in the symplectic case, that the condition {\bf (T4)} 
is satisfied, and the end points are $Q_j=(2n+1,n-j+1+\lambda_j)$ 
since the alphabet is now of length $2n+1$. To take into account the terms in $(1-a_{\lambda_j-j+n+1})$
the lattice paths may now include a final diagonal step.
The fact that it is diagonal ensures that there is at most one of these steps on each
lattice path. Once again it is only the $n$-tuples of non-intersecting lattice paths from $P_i$ to $Q_i$ that
contribute to $so_\lambda(\x,\ovx,1|\a)$ and these are in bijective correspondence with
the odd orthogonal tableaux of Definition~\ref{def-ooTab} of shape $\lambda$ with entries from $\{1<\ov1<\cdots<n<\ov{n}<0\}$.
The fact that on each lattice path the final step is either vertical or diagonal, with the
latter to be associated with entries $0$ ensures that the condition {\bf (T5)} 
is automatically satisfied.
This is exemplified 
for $n=4$ and $\lambda=(4,3,3,0)$  by the following, where we have again used the fact that $a_m=0$ for $m\leq0$ in 
specifying the tableau entry weights but not the lattice path edge weights: 

\begin{equation}\label{eqn-oo-LPTex}
\begin{array}{c}
LP\ =\ 
\vcenter{\hbox{
\begin{tikzpicture}[x={(0in,-0.3in)},y={(0.5in,0in)}] 
\foreach \i in {1,...,4} \draw[-](2*\i-1,5-\i)to(0,4+\i); \foreach \i in {1,...,3} \draw(0-0.2,8-\i+0.2)node{$a_{\ov{\i}}$}; 
\draw(0-0.2,8+0.2)node{$a_{0}$}; 
\foreach \i in {1,...,8} \draw[-](9-0.5,\i-0.5)to(\i-1,9); \foreach \i in {1,...,8} \draw(\i-1.2,9+0.2)node{$a_\i$};
\foreach \i in {1,...,8} \draw[-](9-0.5,\i-0.5)to(9+0.8,\i+0.8); \foreach \i in {1,...,8} \draw(10,\i+0.8)node{$-\!a_\i$};
\draw(1,4-0.4)node{$P_1$};
\draw(3,3-0.4)node{$P_2$};
\draw(5,2-0.4)node{$P_3$};
\draw(7,1-0.4)node{$P_4$};
\draw(9+0.4,1)node{$Q_4$};
\draw(9+0.4,5)node{$Q_3$};
\draw(9+0.4,6)node{$Q_2$};
\draw(9+0.4,8)node{$Q_1$};
\draw[draw=magenta,ultra thick] (1,4)to(1,5); \draw(1-0.3,5-0.5)node{$\magenta{x_1\!\!+\!\!a_{\ov2}}$};
\draw[draw=magenta,ultra thick] (1,5)to(2,5);
\draw[draw=magenta,ultra thick] (2,5)to(2,6); \draw(2-0.3,6-0.5)node{$\magenta{\ov{x}_1\!\!+\!\!a_{0}}$};
\draw[draw=magenta,ultra thick] (2,6)to(3,6);
\draw[draw=magenta,ultra thick] (3,6)to(3,7); \draw(3-0.3,7-0.5)node{$\magenta{x_2\!\!+\!\!a_{2}}$};
\draw[draw=magenta,ultra thick] (3,7)to(8,7);
\draw[draw=magenta,ultra thick] (8,7)to(8,8); \draw(8-0.3,8-0.5)node{$\magenta{\ov{x}_4\!\!+\!\!a_{8}}$};
\draw[draw=magenta,ultra thick] (8,8)to(9,8);
\draw[draw=blue,ultra thick] (3,3)to(5,3);
\draw[draw=blue,ultra thick] (5,3)to(5,4);  \draw(5-0.3,4-0.5)node{$\blue{x_3\!\!+\!\!a_{1}}$};
\draw[draw=blue,ultra thick] (5,4)to(7,4);
\draw[draw=blue,ultra thick] (7,4)to(7,5);  \draw(7-0.3,5-0.5)node{$\blue{x_4\!\!+\!\!a_{4}}$};
\draw[draw=blue,ultra thick] (7,5)to(8,5);  
\draw[draw=blue,ultra thick] (8,5)to(9,6);  \draw(9-0.4,6+0.0)node{$\blue{1\!\!-\!\!a_6}$};
\draw[draw=red,ultra thick] (5,2)to(7,2);
\draw[draw=red,ultra thick] (7,2)to(7,3);  \draw(7-0.3,3-0.5)node{$\red{x_4\!\!+\!\!a_{2}}$};
\draw[draw=red,ultra thick] (7,3)to(8,3);
\draw[draw=red,ultra thick] (8,3)to(8,4);  \draw(8-0.3,4-0.5)node{$\red{\ov{x}_4\!\!+\!\!a_{4}}$};
\draw[draw=red,ultra thick] (8,4)to(9,5);  \draw(9-0.4,5+0.0)node{$\red{1\!\!-\!\!a_5}$};
\draw[draw=cyan,ultra thick] (7,1)to(8,1);
\draw[draw=cyan,ultra thick] (8,1)to(9,1); 
\foreach \j in {4,...,8} \draw(1,\j)node{$\bullet$};
\foreach \j in {4,...,8} \draw(2,\j)node{$\bullet$};
\foreach \j in {3,...,8} \draw(3,\j)node{$\bullet$};
\foreach \j in {3,...,8} \draw(4,\j)node{$\bullet$};
\foreach \j in {2,...,8} \draw(5,\j)node{$\bullet$};
\foreach \j in {2,...,8} \draw(6,\j)node{$\bullet$};
\foreach \j in {1,...,8} \draw(7,\j)node{$\bullet$};
\foreach \j in {1,...,8} \draw(8,\j)node{$\bullet$};
\foreach \j in {1,...,8} \draw(9,\j)node{$\bullet$};
\end{tikzpicture}
}}\cr\cr
T\ =\ \YT{0.2in}{0.2in}{}{
 {\magenta{1},\magenta{\ov1},\magenta{2},\magenta{\ov4}},
 {\blue{3},\blue{4},\blue{0}},
 {\red{4},\red{\ov4},\red{0}},
}
\qquad\qquad
\wgt(T)\ =\ \YT{0.2in}{0.5in}{}{
 {\magenta{x_1},\magenta{\ov{x}_1},\magenta{{x_2}\!+\!{a_2}},\magenta{\ov{x}_4\!+\!{a_8}}},
 {\blue{{x_3}\!+\!{a_{1}}},\blue{{x_4}\!+\!{a_4}},\blue{{1}-{a_6}}},
 {\red{{x_4}\!+\!{a_2}},\red{\ov{x}_4\!+\!{a_4}},\red{{1}-{a_5}}},
} 
\end{array}
\end{equation}

\bigskip
\noindent{\bf Case $o(2n)$}.
This time we choose to base the even orthogonal case on the use of (\ref{eqn-eohm-tau2-n}) which implies that
\begin{equation}\label{eqn-eo-3hm}
  o_\lambda(\x,\ovx|\a) = \left| \begin{array}{c}
              (x_i+a_{i+1-n}) h^{\gl}_{\lambda_j-j+i}(x_i,\x^{(i+1)},\ovx^{(i+1)}\,|\,\tau^{i+1-n}\a) \cr
							+(\ov{x}_i+a_{i+1-n}) h^{\gl}_{\lambda_j-j+i}(\ov{x}_i,\x^{(i+1)},\ovx^{(i+1)}\,|\,\tau^{i+1-n}\a) \cr
							+h^{\gl}_{\lambda_j-j+i}(\x^{(i+1)},\ovx^{(i+1)}\,|\,\tau^{i+1-n}\a) \cr
          \end{array}
					\right|\,.
\end{equation}
This determinant can be expanded as a sum of $3^n$ determinants in each of which one makes a choice of 
which of the three types of entries occur in row $i$ for $i=1,2,\ldots,n$. The non-intersecting lattice path model 
can then be applied to each such determinant independently. In each case the
appropriate underlying lattice is now typically as shown below.
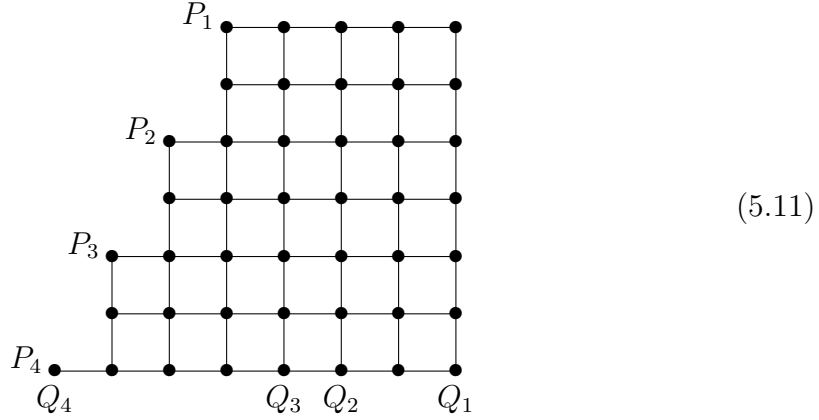
\begin{equation}\label{eqn-eo-grid}
\vcenter{\hbox{
\begin{tikzpicture}[x={(0in,-0.3in)},y={(0.3in,0in)}] 
\foreach \j in {4,...,8} \draw(2,\j)node{$\bullet$};
\foreach \j in {4,...,8} \draw(3,\j)node{$\bullet$};
\foreach \j in {3,...,8} \draw(4,\j)node{$\bullet$};
\foreach \j in {3,...,8} \draw(5,\j)node{$\bullet$};
\foreach \j in {2,...,8} \draw(6,\j)node{$\bullet$};
\foreach \j in {2,...,8} \draw(7,\j)node{$\bullet$};
\foreach \j in {1,...,8} \draw(8,\j)node{$\bullet$};
\draw[-](2,4)to(2,8);
\draw[-](3,4)to(3,8);
\draw[-](4,3)to(4,8);
\draw[-](5,3)to(5,8);
\draw[-](6,2)to(6,8);
\draw[-](7,2)to(7,8);
\draw[-](8,1)to(8,8);
\draw[-](6,2)to(8,2);
\draw[-](4,3)to(8,3);
\foreach \j in {4,...,8} \draw[-](2,\j)to(8,\j);
\draw(2-0.2,3.5)node{$P_1$};
\draw(4-0.2,2.5)node{$P_2$};
\draw(6-0.2,1.5)node{$P_3$};
\draw(8-0.2,0.5)node{$P_4$};
\draw(8.5,1)node{$Q_{4}$};  
\draw(8.5,5)node{$Q_{3}$};
\draw(8.5,6)node{$Q_{2}$};
\draw(8.5,8)node{$Q_{1}$};
\end{tikzpicture}
}}
\end{equation}
This represents a minor modification of the symplectic lattice in which the starting points of the 
lattice paths are still $P_i=(2i-1,n-i+1)$ but the end points are now $Q_j=(2n-1,n-j+1+\lambda_j)$.

However, the major difference is that there are now three distinct types of path and edge labelling corresponding 
to each of the three terms displayed in the determinantal
expression (\ref{eqn-eo-3hm}). For those starting at $P_i=(2i-1,n-i+1)$ the three types are characterised
as follows by their initial edges:\\
\begin{tabular}{rl}
type (i)&a horizontal edge labelled $(x_i+a_{i+1-n})$; \cr 
type (ii)&a horizontal edge labelled $(\ov{x}_i+a_{i+1-n})$; \cr
type (iii)&a vertical edge. 
\end{tabular}
\\
There are further restrictions on succeeding edge weights such that in the case of type (i) $(x_i+a_{i+1-n})$ may be followed by
a sequence $(x_i+a_{i+2-n})(x_i+a_{i+3-n})\cdots$ but no edge of weight $(\ov{x}_i+a_{k})$, whilst in the case of type (ii)
$(\ov{x}_i+a_{i+1-n})$ may be followed by a sequence $(\ov{x}_i+a_{i+2-n})(\ov{x}_i+a_{i+3-n})\cdots$ but no edge of weight $({x}_i+a_{k})$.

A typical non-intersecting $n$-tuple in the case $n=4$ and $\lambda=(5,4,4,4)$ takes the form
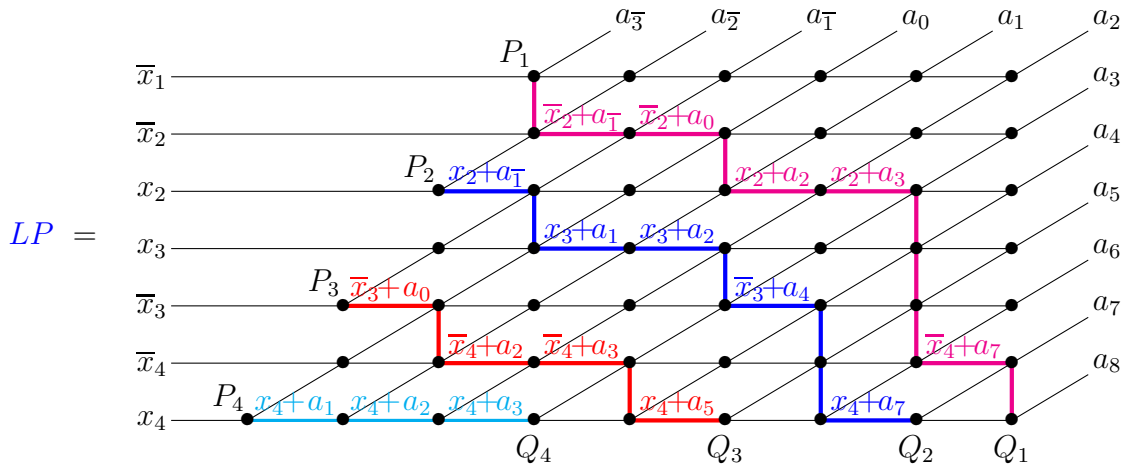
\begin{equation}\label{eqn-eo-LPTex}
\begin{array}{c}
\blue{LP}\ = \ 
\vcenter{\hbox{
\begin{tikzpicture}[x={(0in,-0.3in)},y={(0.5in,0in)}] 
\draw[-](1+0.2,5-0.2)to(2,4); 
\draw[-](1+0.2,6-0.2)to(4,3); 
\draw[-](1+0.2,7-0.2)to(6,2); 
\draw[-](1+0.2,8-0.2)to(8,1); 
\draw[-](1+0.2,9-0.2)to(8,2); 
\foreach \i in {2,...,8} \draw[-,thin](\i-1+0.2,10-0.2)to(8,\i+1); 
\foreach \i in {2,...,8} \draw[-,thin](\i,0.2)to(\i,9);
\draw(2,0)node{$\ov{x}_1$};
\draw(3,0)node{$\ov{x}_2$};
\draw(4,0)node{${x}_2$};
\draw(5,0)node{${x}_3$};
\draw(6,0)node{$\ov{x}_3$};
\draw(7,0)node{$\ov{x}_4$};
\draw(8,0)node{${x}_4$};
\draw(2-0.4,4-0.2)node{$P_1$};
\draw(4-0.4,3-0.2)node{$P_2$};
\draw(6-0.4,2-0.2)node{$P_3$};
\draw(8-0.4,1-0.2)node{$P_4$};
\draw(8+0.5,4)node{$Q_4$};
\draw(8+0.5,6)node{$Q_3$};
\draw(8+0.5,8)node{$Q_2$};
\draw(8+0.5,9)node{$Q_1$};
\draw(1,5)node{$a_{\ov3}$};
\draw(1,6)node{$a_{\ov2}$};
\draw(1,7)node{$a_{\ov1}$};
\draw(1,8)node{$a_0$};
\draw(1,9)node{$a_1$};
\foreach \i in {2,...,8} \draw(\i-1,10)node{$a_\i$};
\draw[draw=magenta,ultra thick] (2,4)to(3,4); 
\draw[draw=magenta,ultra thick] (3,4)to(3,5); \draw(3-0.3,5-0.5)node{$\magenta{\ov{x}_2\!\!+\!\!a_{\ov1}}$};
\draw[draw=magenta,ultra thick] (3,5)to(3,6); \draw(3-0.3,6-0.5)node{$\magenta{\ov{x}_2\!\!+\!\!a_{0}}$};
\draw[draw=magenta,ultra thick] (3,6)to(4,6);
\draw[draw=magenta,ultra thick] (4,6)to(4,7); \draw(4-0.3,7-0.5)node{$\magenta{{x}_2\!\!+\!\!a_2}$};
\draw[draw=magenta,ultra thick] (4,7)to(4,8); \draw(4-0.3,8-0.5)node{$\magenta{{x}_2\!\!+\!\!a_3}$};
\draw[draw=magenta,ultra thick] (4,8)to(7,8); 
\draw[draw=magenta,ultra thick] (7,8)to(7,9); \draw(7-0.3,9-0.5)node{$\magenta{\ov{x}_4\!\!+\!\!a_7}$};
\draw[draw=magenta,ultra thick] (7,9)to(8,9); 
\draw[draw=blue,ultra thick] (4,3)to(4,4);  \draw(4-0.3,4-0.5)node{$\blue{{x}_2\!\!+\!\!a_{\ov1}}$};
\draw[draw=blue,ultra thick] (4,4)to(5,4);
\draw[draw=blue,ultra thick] (5,4)to(5,5);  \draw(5-0.3,5-0.5)node{$\blue{{x}_3\!\!+\!\!a_1}$};
\draw[draw=blue,ultra thick] (5,5)to(5,6);  \draw(5-0.3,6-0.5)node{$\blue{{x}_3\!\!+\!\!a_2}$};
\draw[draw=blue,ultra thick] (5,6)to(6,6);
\draw[draw=blue,ultra thick] (6,6)to(6,7);  \draw(6-0.3,7-0.5)node{$\blue{\ov{x}_3\!\!+\!\!a_4}$};
\draw[draw=blue,ultra thick] (6,7)to(8,7);
\draw[draw=blue,ultra thick] (8,7)to(8,8);  \draw(8-0.3,8-0.5)node{$\blue{{x}_4\!\!+\!\!a_7}$};
\draw[draw=red,ultra thick] (6,2)to(6,3);  \draw(6-0.3,3-0.5)node{$\red{\ov{x}_3\!\!+\!a_{0}}$};
\draw[draw=red,ultra thick] (6,3)to(7,3);  
\draw[draw=red,ultra thick] (7,3)to(7,4);  \draw(7-0.3,4-0.5)node{$\red{\ov{x}_4\!\!+\!\!a_2}$};
\draw[draw=red,ultra thick] (7,4)to(7,5);  \draw(7-0.3,5-0.5)node{$\red{\ov{x}_4\!\!+\!\!a_3}$};
\draw[draw=red,ultra thick] (7,5)to(8,5);  
\draw[draw=red,ultra thick] (8,5)to(8,6);  \draw(8-0.3,6-0.5)node{$\red{{x}_4\!\!+\!\!a_5}$};
\draw[draw=cyan,very thick] (8,1)to(8,2);  \draw(8-0.3,2-0.5)node{$\cyan{x_4\!\!+\!a_1}$};
\draw[draw=cyan,very thick] (8,2)to(8,3);  \draw(8-0.3,3-0.5)node{$\cyan{{x}_4\!\!+\!a_2}$};
\draw[draw=cyan,very thick] (8,3)to(8,4);  \draw(8-0.3,4-0.5)node{$\cyan{{x}_4\!\!+\!\!a_3}$};
\foreach \j in {4,...,9} \draw(2,\j)node{$\bullet$};
\foreach \j in {4,...,9} \draw(3,\j)node{$\bullet$};
\foreach \j in {3,...,9} \draw(4,\j)node{$\bullet$};
\foreach \j in {3,...,9} \draw(5,\j)node{$\bullet$};
\foreach \j in {2,...,9} \draw(6,\j)node{$\bullet$};
\foreach \j in {2,...,9} \draw(7,\j)node{$\bullet$};
\foreach \j in {1,...,9} \draw(8,\j)node{$\bullet$};
\end{tikzpicture}
}}\cr\cr
\end{array}
\end{equation}

In this example $P_1Q_1$ is of type (iii), while the paths $P_2Q_2$, $P_3Q_3$ and $P_4Q_4$ 
have been chosen to be of types (i), (ii) and (i), respectively. The paths represent contributions to 
the diagonal elements $h_5^{\gl}(\ov{x}_2,{x}_2,{x}_3,\ov{x}_3,\ov{x}_4,{x}_4|\tau^{-2}\a)$,
$(x_2+a_{\ov1})h_4^{\gl}(x_2,{x}_3,\ov{x}_3,\ov{x}_4,{x}_4|\tau^{-1}\a)$, $(\ov{x}_3+a_{0})h_{3}^{\gl}(\ov{x}_3,\ov{x}_4,{x}_4|\a)$
and $(x_4+a_{1})h_{2}^{\gl}({x}_4|\tau\a)$, where advantage has been taken of the symmetry of each $h_m^{\gl}$
with respect to its arguments in order to write these arguments in the order dictated by our choice of paths of particular types.
It is this choice that has forced the first 
edge labels at levels $3$, $5$ and $7$, counted from top to bottom, to involve ${x}_2$, $\ov{x}_3$ and $x_4$,
respectively. This labelling is maintained across each fixed level, with labels involving their
inverses $\ov{x}_2$, ${x}_3$ and $\ov{x}_4$ applied to all edges at levels $2$, $4$ and $6$. 
Thanks to the path $P_1Q_1$ being of type (iii), the edge labelling at level $1$ is irrelevant in this example
but has been chosen to be $\ov{x}_1$.

It follows from the usual non-intersecting paths argument in respect of each determinant, that in the even orthogonal
case $o(2n)$ the character $o_\lambda(\x,\ovx|\a)$ may be evaluated by enumerating all possible non-intersecting $n$-tuples
of lattice paths between $P_i=(2i-1;n-i+1)$ and $Q_j = (2n-1;n-j+1+\lambda_j)$ on the even orthogonal lattice of the 
type illustrated in (\ref{eqn-eo-grid}). The weight $(z_k+a_l)$ attached to each horizontal edge is determined 
in the usual way by the labels $z_k$ and $a_l$, respectively, on the horizontal and diagonal lines intersecting at its right hand end.
The default horizontal labels are $\ov{x}_i$ at the level of $P_i$ for $i=1,2,\ldots,n$ and $x_i$ at the
level immediately above $P_i$ for $i=2,3,\ldots,n$. 
The diagonal through $P_n$ is labelled $a_0$, with consecutive diagonals to the right and left labelled $a_1,a_2,\ldots$ and 
$a_{\ov1},a_{\ov2},\ldots$, respectively, although in the final evaluation of $o_\lambda(\x,\ovx|\a)$ one must set
$a_\ell=0$ for all $\ell\leq0$.
However, having evaluated the contribution of such an $n$-tuple in this way one must then carry out a summation over all possible independent
interchanges of horizontal labels $x_k$ and $\ov{x}_k$ for those $k$ for which the first edge of the path starting from $P_k$
is horizontal.

One result of this freedom to interchange horizontal labels is that
the normal passage from a lattice path $n$-tuple to a tableau, obtained by reading consecutive labels along each path
$P_iQ_i$ and entering them successively from left to right across the $i$th row, yields
a tableau $\tilde{T}$ that does not necessarily belong to ${\cal T}_\lambda^{\eo}$ 
since it may involve an alphabet of entries in non-standard order. 
In the case of our example (\ref{eqn-eo-LPTex}) the non-standard order is $1<\ov1<\ov2<2<3<\ov3<\ov4<4$, 
as determined by the horizontal labels, and the resulting tableau $\tilde{T}$ and its weighting are as shown below:
\begin{equation}\label{eqn-eo-tildeT}
\blue{\widetilde{T}}\ = \
\YT{0.2in}{0.2in}{}{
 {\magenta{\ov2},\magenta{\ov2},\magenta{2},\magenta{2},\magenta{\ov4}},
 {\blue{2},\blue{3},\blue{3},\blue{\ov3},\blue{4}},
 {\red{\ov3},\red{\ov4},\red{\ov4},\red{4}},
 {\cyan{4},\cyan{4},\cyan{4}},
}
\qquad\qquad
\wgt(\blue{\widetilde{T}})\ = \ 
\YT{0.2in}{0.5in}{}{
 {\magenta{\ov{x}_2\!+\!a_{\ov1}},\magenta{\ov{x}_2\!+\!a_0},\magenta{{x}_2\!+\!a_2},\magenta{{x}_2\!+\!a_3},\magenta{\ov{x}_{4}\!+\!{a_7}}},
 {\blue{{x}_2\!+\!a_{\ov1}},\blue{{x}_3\!+\!{a_1}},\blue{{x}_3\!+\!{a_2}},\blue{\ov{x}_3\!+\!a_4},\blue{{x}_4\!+\!a_7}},
 {\red{\ov{x}_3\!+\!a_0},\red{\ov{x}_{4}\!+\!{a_2}},\red{\ov{x}_4\!+\!{a_3}},\red{{x}_4\!+\!{a_5}}},
 {\cyan{{x_4}\!+\!{a_1}},\cyan{x_{4}\!+\!{a_2}},\cyan{x_{4}\!+\!{a_3}}},
} 
\end{equation}

To remedy this situation in which $\widetilde{T}\notin{\cal T}_\lambda^{\eo}$ we move to a standard ordering
of entries by applying an operator $\psi$ to $\widetilde{T}$, where $\psi=\prod_k \chi(\widetilde{T}_{k1}=k)\,\psi_k$
and the action of $\psi_k$ is to interchange all entries $k$ and $\ov{k}$ except for those entries $k$ in the
$k$th row that are not covered by an entry $\ov{k}$ immediately above it. It might be noted that the truth function 
$\chi(\widetilde{T}_{k1}=k)=1$ if and only if the corresponding $n$-tuple of lattice paths $LP$ from which $\widetilde{T}$ 
is derived contains a path $P_kQ_k$ of type (i). In our example this means that $\psi=\psi_2\,\psi_4$, and applying 
this to $\widetilde{T}$ yields $T$ as shown below on the left:
\begin{equation}\label{eqn-eo-TwgtT}
T\ = \
\YT{0.2in}{0.2in}{}{
 {\magenta{2},\magenta{2},\magenta{\ov2},\magenta{\ov2},\magenta{4}},
 {\blue{\ov2},\blue{3},\blue{3},\blue{\ov3},\blue{\ov4}},
 {\red{\ov3},\red{4},\red{4},\red{\ov4}},
 {\cyan{4},\cyan{\ov4},\cyan{\ov4}},
}
\qquad\qquad
\wgt(T)\ = \ 
\YT{0.2in}{0.5in}{}{
 {\magenta{{x}_2\!+\!a_{\ov1}},\magenta{{x}_2\!+\!a_0},\magenta{\ov{x}_2\!+\!a_2},\magenta{\ov{x}_2\!+\!a_3},\magenta{{x}_{4}\!+\!{a_7}}},
 {\blue{\ov{x}_2\!+\!a_{\ov1}},\blue{{x}_3\!+\!{a_1}},\blue{{x}_3\!+\!{a_2}},\blue{\ov{x}_3\!+\!a_4},\blue{\ov{x}_4\!+\!a_7}},
 {\red{\ov{x}_3\!+\!a_0},\red{{x}_{4}\!+\!{a_2}},\red{{x}_4\!+\!{a_3}},\red{\ov{x}_4\!+\!{a_5}}},
 {\cyan{{x_4}\!+\!{a_1}},\cyan{\ov{x}_{4}\!+\!{a_2}},\cyan{\ov{x}_{4}\!+\!{a_3}}},
} 
\end{equation}
As intended, this transformation $\psi$ of $\widetilde{T}$ leads to a tableau $T\in{\cal T}_\lambda^{\eo}$.
In particular, the entries are now weakly increasing across rows and strictly increasing down columns with respect 
to the standard order $1<\ov1<2<\ov2<3<\ov3<4<\ov4$, with both $k$ and $\ov{k}$ appearing no lower than the
$k$th row, in accordance with {\bf (T4)}, while the condition {\bf (T6)} is automatically satisfied since 
the only way in which entries $\ov{k}$ can appear in the $k$th row of $T$ to the right of an entry $k$ is 
if each such entry was originally immediately above an entry $k$ in $\widetilde{T}$ with which it was exchanged 
under the action of $\psi_k$.

On the right of (\ref{eqn-eo-TwgtT}) we have shown $\wgt(T)$ which is obtained from $\wgt(\widetilde{T})$ by
simply interchanging $x_k$ and $\ov{x}_k$ wherever $k$ and $\ov{k}$ have been interchanged in 
passing from $\widetilde{T}$ to $T=\psi\,\widetilde{T}$. Unfortunately, 
this type of interchange is not in general weight preserving, as can be seen in this example. However, the transformation $\psi$ is weight 
preserving when applied to certain sums of tableaux. In our case the required combination of non-standard 
tableaux $\widetilde{T}$ is 
\begin{equation}
\begin{array}{cccc}
\phantom{+}\YT{0.2in}{0.2in}{}{
 {\magenta{\ov2},\magenta{\ov2},\magenta{\ov2},\magenta{\ov2},\magenta{\ov4}},
 {\blue{2},\blue{3},\blue{3},\blue{\ov3},\blue{4}},
 {\red{\ov3},\red{4},\red{4},\red{\ov4}},
 {\cyan{4},\cyan{4},\cyan{4}},
}
&+
\YT{0.2in}{0.2in}{}{
 {\magenta{\ov2},\magenta{\ov2},\magenta{\ov2},\magenta{2},\magenta{\ov4}},
 {\blue{2},\blue{3},\blue{3},\blue{\ov3},\blue{4}},
 {\red{\ov3},\red{\ov4},\red{\ov4},\red{\ov4}},
 {\cyan{4},\cyan{4},\cyan{4}},
}
&+
\YT{0.2in}{0.2in}{}{
 {\magenta{\ov2},\magenta{\ov2},\magenta{2},\magenta{2},\magenta{\ov4}},
 {\blue{2},\blue{3},\blue{3},\blue{\ov3},\blue{4}},
 {\red{\ov3},\red{\ov4},\red{\ov4},\red{\ov4}},
 {\cyan{4},\cyan{4},\cyan{4}},
}
&+
\YT{0.2in}{0.2in}{}{
 {\magenta{\ov2},\magenta{2},\magenta{2},\magenta{2},\magenta{\ov4}},
 {\blue{2},\blue{3},\blue{3},\blue{\ov3},\blue{4}},
 {\red{\ov3},\red{\ov4},\red{\ov4},\red{\ov4}},
 {\cyan{4},\cyan{4},\cyan{4}},
}
\cr\cr
+\YT{0.2in}{0.2in}{}{
 {\magenta{\ov2},\magenta{\ov2},\magenta{\ov2},\magenta{\ov2},\magenta{\ov4}},
 {\blue{2},\blue{3},\blue{3},\blue{\ov3},\blue{4}},
 {\red{\ov3},\red{\ov4},\red{\ov4},\red{4}},
 {\cyan{4},\cyan{4},\cyan{4}},
}
&+
\YT{0.2in}{0.2in}{}{
 {\magenta{\ov2},\magenta{\ov2},\magenta{\ov2},\magenta{2},\magenta{\ov4}},
 {\blue{2},\blue{3},\blue{3},\blue{\ov3},\blue{4}},
 {\red{\ov3},\red{\ov4},\red{\ov4},\red{4}},
 {\cyan{4},\cyan{4},\cyan{4}},
}
&+
\YT{0.2in}{0.2in}{}{
 {\magenta{\ov2},\magenta{\ov2},\magenta{2},\magenta{2},\magenta{\ov4}},
 {\blue{2},\blue{3},\blue{3},\blue{\ov3},\blue{4}},
 {\red{\ov3},\red{\ov4},\red{\ov4},\red{4}},
 {\cyan{4},\cyan{4},\cyan{4}},
}
&+
\YT{0.2in}{0.2in}{}{
 {\magenta{\ov2},\magenta{2},\magenta{2},\magenta{2},\magenta{\ov4}},
 {\blue{2},\blue{3},\blue{3},\blue{\ov3},\blue{4}},
 {\red{\ov3},\red{\ov4},\red{\ov4},\red{4}},
 {\cyan{4},\cyan{4},\cyan{4}},
}
\cr
\end{array}
\end{equation}
where in each of the two rows of this sum triples of entries $\ov2\ov2\ov2$, $\ov2\ov22$, $\ov222$ or $222$ appear
in the first row of each tableau to the right of the vertical pair $\ov2$ over $2$, while the two rows are distinguished only 
by the single entries $4$ and $\ov4$ at the end of the third row of each tableau immediately to the right of a vertical pair $\ov4$ over $4$.
The action of $\psi=\psi_2\psi_4$ inverts all the vertical pairs $\ov2$ over $2$ and $\ov4$ over $4$, replaces the above triples
by $222$, $22\ov2$, $2\ov2\ov2$ or $\ov2\ov2\ov2$, respectively, and interchanges the single entries $4$ and $\ov4$.

As far as the corresponding weights are concerned, each vertical pair $\ov{k}$ over $k$ for $k=2$ and $k=4$, carries a combined weight of 
$(x_k+a_\ell)(\ov{x}_k+a_\ell)$ for some $\ell$. This is obviously invariant under the interchange of $\ov{x}_k$ and $x_k$. 
Furthermore, the sum of the contributions to the weight of the triples of $\ov2$s and $2$s yields a factor
\begin{align}
        &~~(\ov{x}_2+a_0)(\ov{x}_2+a_1)(\ov{x}_2+a_2)+(\ov{x}_2+a_0)(\ov{x}_2+a_1)({x}_2+a_3)\cr
           &+ (\ov{x}_2+a_0)({x}_2+a_2)({x}_2+a_3)+({x}_2+a_1)({x}_2+a_2)({x}_2+a_3)\,.        
\end{align}
It is not immediately obvious that this is invariant under the interchange of $x_2$ and $\ov{x}_2$, but this is the case
since it can be identified with $h_3^{\gl}(\ov{x}_2,x_2|\tau^{-1}\a)=h_3^{\gl}({x}_2,\ov{x}_2|\tau^{-1}\a)$.
Finally, the sum of the contributions from the single entries $4$ and $\ov4$ gives a factor $(\ov{x}_4+a_4)+(x_4+a_5)$
which is again clearly invariant under the interchange of $x_4$ and $\ov{x}_4$. It might also be noted that
this factor can be dentified with $h_1^{\gl}(\ov{x}_4,x_4|\tau^3\a)=h_1^{\gl}(x_4,\ov{x}_4|\tau^3\a)$.

More generally, this identification of common factors with $h_m^{\gl}(x_k,\ov{x}_k|\tau^r\a)$ for some $m,k,r$ 
is what allows us to conclude that the action of $\psi_k$ is weight preserving on sums of non-standard tableaux. 
This comes about because in any tableau, $\widetilde{T}$, standard with respect to an order $\cdots<\ov{k}<k<\cdots$
the entries $k$ and $\ov{k}$ are distributed across any row in a sequence of say $p$ entries $\ov{k}$
followed by $m-p$ entries $k$ that may or may not be sandwiched between any numbers of $\ov{k}$ over $k$ pairs at either end,
as shown schematically below
\begin{equation}
\vcenter{\hbox{
\begin{tikzpicture}[x={(0in,-0.3in)},y={(0.3in,0in)}] 
\foreach \i/\j in {{0/8},{1/1},{1/2},{1/4}} \draw (\i,\j) rectangle +(-1,-1); 
\foreach \i/\j in {{0/8},{1/1},{1/2},{1/4}} \draw (\i-0.5,\j-0.5) node {$\ov{k}$};
\foreach \i/\j in {{2/1},{1/5},{1/7},{1/8}} \draw (\i,\j) rectangle +(-1,-1);  
\foreach \i/\j in {{2/1},{1/5},{1/7},{1/8}} \draw (\i-0.5,\j-0.5) node {$k$};
\foreach \i/\j in {{1/3},{1/6}} \draw (\i,\j) rectangle +(-1,-1);  
\foreach \i/\j in {{1/3},{1/6}} \draw (\i-0.5,\j-0.5) node {$\cdots$};
\draw[latex-](1.5,1)--(1.5,3-0.7); \draw(1.5,3-0.5)node{$p$}; \draw[-latex](1.5,3-0.2)--(1.5,4); 
\draw[latex-](1.5,4)--(1.5,5-0.2); \draw(1.5,6-0.5)node{$m\!-\!p$}; \draw[-latex](1.5,6+0.2)--(1.5,7);

\end{tikzpicture}
}}
\end{equation}
The corresponding weights of the entries in such a portion of $\widetilde{T}$ take the form 
\begin{equation}
\vcenter{\hbox{
\begin{tikzpicture}[x={(0in,-0.3in)},y={(0.42in,0in)}] 
\foreach \i/\j in {{0/13},{1/1},{1/3},{1/6},{1/8},{1/11},{1/13},{2/1}} \draw (\i,\j) rectangle +(1,2); 
\draw (1+0.5,1+1) node {$\ov{x}_k\!+\!a_r$};\draw (2+0.5,1+1) node {$x_k\!+\!a_r$};
\draw (1+0.5,3+1) node {$\ov{x}_k\!+\!a_{r+1}$};\draw (1+0.5,6+1) node {$\ov{x}_k\!+\!a_{r+p}$};
\draw (1+0.5,8+1) node {${x}_k\!+\!a_{r+p+2}$};\draw (1+0.5,11+1) node {${x}_k\!+\!a_{r+m+1}$};
\draw (0+0.5,13+1) node {$\ov{x}_k\!+\!a_{r+m+2}$};\draw (1+0.5,13+1) node {${x}_k\!+\!a_{r+m+2}$};
\foreach \i/\j in {{1/5},{1/10}} \draw (\i,\j) rectangle +(1,1); 
\foreach \i/\j in {{1/5},{1/10}} \draw (\i+0.5,\j+0.5) node {$\cdots$};
\draw[latex-](2.5,3)--(2.5,5+0.2); \draw(2.5,5.5)node{$p$}; \draw[-latex](2.5,6-0.2)--(2.5,8);
\draw[latex-](2.5,8)--(2.5,10-0.1); \draw(2.5,10.5)node{$m\!-\!p$}; \draw[-latex](2.5,11+0.1)--(2.5,13);
\end{tikzpicture}
}}
\end{equation}
The product $(\ov{x}_k+a_l)(x_k+a_l)$ of weights arising from each vertical pair is symmetric in $x_k$ and $\ov{x}_k$,
but that arising from the horizontal sequence of length $m$ is not symmetric. However, it may be recognised as taking the
the form of a single monomial summand of $h_m^{\gl}(\ov{x}_k,x_k|\tau^s\a)$, and summing over $p$ from $0$ to $m$ gives a weight contribution 
precisely equal to $h_m^{\gl}(\ov{x}_k,x_k|\tau^s\a)$. This expression is symmetric under the interchange of $x_k$ and $\ov{x}_k$, so  
that the interchange of $\ov{x}_k$ and $x_k$ induced by the action of $\psi_k$ is weight preserving when acting on 
the corresponding sum of $\widetilde{T}$.

It will be noted that the $T\in{\cal T}_\lambda^{\eo}$ obtained in this way through the action $\psi_k$
on some $\widetilde{T}$ with $\ov{k}$ immediately above an entry $k$ in position $(k,1)$ is indistinguishable 
from some other $\widetilde{T}'$ with $k$ immediately above an entry $\ov{k}$ in position $(k,1)$ to which it is not necessary to
apply $\psi_k$. This results in a doubling of the weight of a given tableau in ${\cal T}_\lambda^{eo}$ for each pair of entries 
$k$ immediately above $\ov{k}$ in the first column of rows $k-1$ and $k$, respectively. 

This completes the proof of the tableau expression for $o_\lambda(\x,\ovx|\a)$ tabulated in Theorem~\ref{thm-Twgts}. 
\qed

\section{Factorial irreducible characters of $SO(2n,\C)$}
\label{sec-so2n}

The irreducible representation of $O(2n,\C)$ of highest weight $\lambda$ remains irrreducible on restriction 
to $SO(2n,\C)$ if $\ell(\lambda)<n$, but reduces to a sum of a pair of irreducible representations of $SO(2n,\C)$
of highest weights $\lambda_+$ and $\lambda_-$ if $\ell(\lambda)=n$, where $\lambda_+=(\lambda_1,\ldots,\lambda_{n-1},\lambda_n)$
and $\lambda_-=(\lambda_1,\ldots,\lambda_{n-1},-\lambda_n)$. 
In terms of characters we have
\begin{align}
o_\lambda(\x,\ovx)&= so_{\lambda}(\x,\ovx)  
                          \quad\hbox{if $\ell(\lambda)<n$}\,; \label{eqn-o-so} \\
o_\lambda(\x,\ovx)&= so_{\lambda_{+}}(\x,\ovx)+so_{\lambda_{-}}(\x,\ovx)  
                          \quad\hbox{if $\ell(\lambda)=n$}\,; \label{eqn-o-so-char}
\end{align}
where~\cite{Lit50,FH04}
\begin{align}
so_\lambda(\x,\ovx)
&=  \frac{\frac12\left|\, x_i^{\lambda_j+n-j} + \ov{x}_i^{\lambda_j+n-j} \,\right| }
		                          {\frac12\left|\, x_i^{n-j} + \ov{x}_i^{n-j} \,\right|} \quad\mbox{if $\ell(\lambda)<n$}\,; \label{eqn-seo} \\
so_{\lambda_{\pm}}(\x,\ovx)
&=  \frac{\frac12\left|\, x_i^{\lambda_j+n-j} + \ov{x}_i^{\lambda_j+n-j} \,\right| \,\pm \, \frac12\left|\, x_i^{\lambda_j+n-j} - \ov{x}_i^{\lambda_j+n-j} \,\right|}
		                          {\frac12\left|\, x_i^{n-j} + \ov{x}_i^{n-j} \,\right|} \quad\mbox{if $\ell(\lambda)=n$}\,. \label{eqn-seo-char} 
\end{align}

In the case $\ell(\lambda)=n$ we introduce the so-called difference characters of Littlewood~\cite{Lit50}:
\begin{equation}
o'_\lambda(\x,\ovx)= so_{\lambda_{+}}(\x,\ovx)-so_{\lambda_{-}}(\x,\ovx)\,. \label{eqn-seo-diffchar} 
\end{equation}
In terms of these we have
\begin{equation}
   so_{\lambda_\pm} = \frac12 \left(\, o_\lambda(\x,\ovx) \pm o'_\lambda(\x,\ovx)\, \right) \label{eqn-sopm}
\end{equation}
with
\begin{equation}
o_\lambda(\x,\ovx) =  \frac{ \left|\, x_i^{\lambda_j+n-j} + \ov{x}_i^{\lambda_j+n-j} \,\right|}
		                          {\frac12\left|\, x_i^{n-j} + \ov{x}_i^{n-j} \,\right|}
\quad\mbox{and}\quad
o'_\lambda(\x,\ovx) =  \frac{ \left|\, x_i^{\lambda_j+n-j} - \ov{x}_i^{\lambda_j+n-j} \,\right|}
                           {\frac12\left|\, x_i^{n-j} + \ov{x}_i^{n-j} \,\right|} \,.   \label{eqn-seo-oo'}
\end{equation}

The passage to factorial versions of these characters and difference characters is accomplished through the 
the following:
\begin{Definition} \label{def-soxxa}
\begin{equation}
so_{\lambda_{\pm}}(\x,\ovx|\a)
=  \frac12 \left(\, o_{\lambda}(\x,\ovx|\a) \pm o'_{\lambda}(\x,\ovx|\a)\,\right)\,, \label{eqn-seoa} 
\end{equation}
where 
\begin{align}
o_{\lambda}(\x,\ovx|\a)
&= \frac{ \left|\, (x_i|\a)^{\lambda_j+n-j} + (\ov{x}_i|\a)^{\lambda_j+n-j} \,\right|}
		                        {\frac12\left|\, (x_i|\a)^{n-j} + (\ov{x}_i|\a)^{n-j} \,\right|}\,; \label{eqn-oxxa} \\
o'_{\lambda}(\x,\ovx|\a)
&=  \frac{ \left|\, (x_i|\a)^{\lambda_j+n-j} - (\ov{x}_i|\a)^{\lambda_j+n-j} \,\right|}
		                          {\frac12\left|\, (x_i|\a)^{n-j} + (\ov{x}_i|\a)^{n-j} \,\right|}\,. \label{eqn-o'xxa}
\end{align}
\end{Definition} 

In order to establish a combinatorial model for factorial difference characters we first show that they satisfy 
a flagged Jacobi-Trudi identity. To this end we first make the following
\begin{Definition}
For any integer $m$ let
\begin{equation}
   h_m^{\eod}(\x,\ovx|\a) = [t^m]\ \left( \frac{1}{1-tx_1}-\frac{1}{1-t\ov{x}_1}\right)\, \prod_{i=2}^n \frac{1}{(1-tx_i)(1-t\ov{x}_i)} \prod_{j=1}^{m+n-1} (1+ta_j)\,.
	    \label{eqn-eod-hma}
\end{equation}
Then for all $m\leq 0$ we have $h_m^{\eod}(\x,\ovx|\a)=0$. 
\end{Definition}

In the case $n=1$ and $\x=(x_i)$ this gives
\begin{equation}
  h_m^{\eod}(x_i,\ov{x}_i|\a) =[t^m]\ \left( \frac{1}{1-tx_i}-\frac{1}{1-t\ov{x}_i}\right)\, \prod_{j=1}^{m} (1+ta_j)=(x_i|\a)^m-(\ov{x}_i|\a)^m\,, \label{eqn-xxbar-hm}
\end{equation}
so that
\begin{equation}\label{eqn-o'det-heod}
    o'_{\lambda}(\x,\ovx|\a)= \prod_{1\leq i<j\leq n}\frac{1}{x_i+\ov{x}_i-x_j-\ov{x}} \
		                      \left|\,  h^{\eod}_{\lambda_j+n-j}(x_i,\ov{x}_i|\a) \,\right|\,,	                          
\end{equation}
where use has been made of the $\a$ independent denominator identity (\ref{eqn-eo-denom}).
However 
\begin{equation}\label{eqn-txxbar}
\frac{1}{1-tx_i}-\frac{1}{1-t\ov{x}_i}=\frac{tx_i}{1-tx_i}-\frac{t\ov{x}_i}{1-t\ov{x}_i}=\frac{t(x_i-\ov{x}_i)}{(1-tx_i)(1-t\ov{x}_i)}\,.
\end{equation}
Extracting common factors of $(x_i-\ov{x}_i)$ from the determinant in (\ref{eqn-o'det-heod}) gives
\begin{equation}
    o'_{\lambda}(\x,\ovx|\a)= \frac{\prod_{i=1}^n (x_i-\ov{x}_i)}{\prod_{1\leq i<j\leq n} (x_i+\ov{x}_i-x_j-\ov{x})}\
		                      \left|\,  h^{\spd}_{\lambda_j+n-j}(x_i,\ov{x}_i|\a) \,\right|\,, \label{eqn-o'-hspd}                          
\end{equation}
where
\begin{equation}
h^{\spd}_m(\x,\ovx|\a) = [t^m]\ t\,\prod_{i=1}^n \frac{1}{(1-tx_i)(1-t\ov{x}_i)}\  \prod_{j=1}^{m+n-1}(1+ta_j) \label{eqn-eohspd}
\end{equation}
Since the only difference between $h^{\spd}_m(\x,\ovx|\a)$ and $h^{\sp}_m(\x,\ovx|\a)$ is an extra factor of $t$, they both satisfy the same recurrence relation
(\ref{eqn-sp-hra-xixj}). This can then be used as in (\ref{eqn-sp-hlambda}) to show that
\begin{equation}
  \left|\,h^{\spd}_{\lambda_j+n-j}(x_i,\ov{x}_i|\a)\,\right|=\prod_{1\leq i<j\leq n} (x_i+\ov{x}_i-x_j-\ov{x}_j) \ 
	                        \left|\,h^{\spd}_{\lambda_j+n-j}(\x^{(i)},\ovx^{(i)}|\a)\,\right|\,.
\end{equation}
Using this in (\ref{eqn-o'-hspd}) and reinserting the factors $(x_i-\ov{x}_i)$ gives 
\begin{align}
o'_\lambda(\x,\ovx|\a) &= \prod_{i=1}^n (x_i-\ov{x}_i) \ \left|\,h^{\spd}_{\lambda_j+n-j}(\x^{(i)},\ovx^{(i)}|\a)\,\right| 
          = \left|\,(x_i-\ov{x}_i) \ h^{\spd}_{\lambda_j+n-j}(\x^{(i)},\ovx^{(i)}|\a)\,\right|\cr 
          &= \left|\,\ h^{\eod}_{\lambda_j+n-j}(\x^{(i)},\ovx^{(i)}|\a)\,\right|\,. \label{eqn-fJT-o'}
\end{align}
This is our required flagged Jacobi-Trudi identity for factorial difference characters.

Then, following the same procedure as in (\ref{eqn-eohm-tau2-n}), we may use (\ref{eqn-txxbar}) to see that for $n>1$
\begin{align}
&h^{\eod}_{m}(\x,\ovx|\a) 
    = [t^m]\ \left(\frac{t(x_1-\ov{x}_1}{(1-tx_1)(1-t\ov{x}_1}\right)\ \prod_{i=2}^n \frac{1}{(1-tx_i)(1-t\ov{x}_i)}\  \prod_{j=1}^{m+n-1}(1+ta_j)\cr
	 & = x_1\,[t^{m-1}]\ \frac{1}{1-tx_1}\ \prod_{i=2}^n \frac{1}{(1-tx_i)(1-t\ov{x}_i)}\  \prod_{j=1+(2-n)}^{(m-1+2n-1-1)+(2-n)}(1+ta_j)\cr 
	 & + \ov{x}_1\,[t^{m-1}]\ \frac{1}{1-t\ov{x}_1}\ \prod_{i=2}^n \frac{1}{(1-tx_i)(1-t\ov{x}_i)}\  \prod_{j=1+(2-n)}^{(m-1+2n-1-1)+(2-n)}(1+ta_j)\cr 
	 &=~~~~ (x_1+a_{2-n}) h^{\gl}_{m-1}(x_1,x_2,\ov{x}_2,\ldots,x_n,\ov{x}_n\,|\,\tau^{2-n}\a)\cr\cr
	 &~~~~ - (\ov{x}_1+a_{2-n}) h^{\gl}_{m-1}(\ov{x}_1,x_2,\ov{x}_2,\ldots,x_n,\ov{x}_n\,|\,\tau^{2-n}\a)\,, \label{eqn-eohmd-tau2-n}
\end{align} 
where we have used the fact that $a_{2-n}=0$ for $n>1$. The same result applies in the case $n=1$ as can be seen
from (\ref{eqn-xxbar-hm}). 

This enables us to set up combinatorial models for $o'_\lambda(\x,\ovx|\a)$ first in terms of non-intersecting lattice paths and then 
tableaux exactly as was done for $o_\lambda(\x,\ovx|\a)$. This leads to
\begin{Theorem}\label{thm-Twgts-eod}
 Let $\lambda$ be a partition of length $\ell(\lambda)=n$ and let ${\cal T}^{\eod}_\lambda$ be the subset of 
${\cal T}^{\eod}_\lambda$ such that $T_{k1}\in\{k,\ov{k}\}$ for all $k=1,2,\ldots,n$ and let $\bar(T)$ be the number of
barred entries in the first column of $T$. Then
\begin{equation}\label{eqn-diffchar-tab}
   o'_\lambda(\x,\ovx|\a) = \sum_{T\in{\cal T}^{\eod}_\lambda}\ (-1)^{\bar(T)} \ \prod_{(i,j)\in F^\lambda} \wgt(T_{ij})\,
\end{equation} 
where
\begin{equation}\label{eqn-Twgts-eodn}
\begin{array}{|l|l|l|l|}
\hline
T_{ij}&\wgt(T_{ij})&\cr
\hline
k&x_k+a_{2k-1+n-j-i+\delta_{ik}}&a_m=0\hbox{~~for~~}m\leq0\cr
 \ov{k}&\ov{x}_k+a_{2k+n-j-i}&\cr		
\hline	
\end{array}
\end{equation}																														
\end{Theorem}
\medskip

\noindent{\bf Proof}:
Since $\ell(\lambda)=n$ we are led by virtue of (\ref{eqn-eohmd-tau2-n}) to only those $n$-tuples 
of non-intersecting lattice paths $P_iQ_i$ with $i=1,2,\ldots,n$ involving just two types of lattice path, namely those previously referred to 
as type (i) or type (ii). These necessarily begin in the case of $P_kQ_k$ with horizontal edges labelled $(x_k+a_{k+1-n})$ or $(\ov{x}_k+a_{k+1-n})$, respectively.
After standardisation through the action of $\psi$, as described in dealing with $o_\lambda(\x,\ovx|\a)$, this restricts the allowed 
tableaux to those in the subset ${\cal T}^{\eod}_\lambda$, as required. 
Noting that the number, $\zeta(T)$, of pairs of entries $T_{k-1,1}=k$ and $T_{k1}=\ov{k}$
in $T\in{\cal T}^{\eod}_\lambda$ is necessarily $0$, and taking into account the sign of the second term in (\ref{eqn-eohmd-tau2-n}),
one is led immediately to (\ref{eqn-diffchar-tab}), where the values of the weights in the last row of
(\ref{eqn-Twgts}) remain the same.
\qed

Combining the results of Theorem~\ref{thm-Twgts} in the case of $o_\lambda(\x,\ovx|\a)$ and the above Theorem~\ref{thm-Twgts-eod} for
$o'_\lambda(\x,\ovx|\a)$ one arrives at the following theorem covering all factorial irreducible characters of $so(n)$ of highest weight
$\lambda$ with $\ell(\lambda)<n$ or $\lambda_\pm$ with $\ell(\lambda)=n$

\begin{Theorem}\label{thm-Twgts-so2n}
For each $T\in{\cal T}^{\eo}_\lambda$ let $\bar(T)$ be the number of entries $T_{k,1}=\ov{k}$ in the first column of $T$,
and let $\zeta(T)$ be the number of pairs of entries $T_{k-1,1}=k$ and $T_{k,1}=\ov{k}$ appearing in
the first column of $T$. Then for $\ell(\lambda)<n$ 
\begin{equation}\label{eqn-wgt-tab-o2n}
   so_\lambda(\x,\ovx|\a) = \sum_{T\in{\cal T}^{\eo}_\lambda}\ 2^{\zeta(T)}\ \prod_{(i,j)\in F^\lambda} \wgt(T_{ij})\,
\end{equation}
and for $\ell(\lambda)=n$
\begin{equation}\label{eqn-wgt-tab-so2n}
   so_{\lambda_\pm}(\x,\ovx|\a) = \sum_{T\in{\cal T}^{\eo}_\lambda}\ \left(\, \theta_{\zeta(T)}\ 2^{\zeta(T)-1}
	                                 + \delta_{\zeta(T),1}\ \epsilon_{\pm}(T)\,\right)\ \prod_{(i,j)\in F^\lambda} \wgt(T_{ij})\,.
\end{equation}
where $\theta_{\zeta(T)}=1$ if $\zeta(T)>0$ and $0$ otherwise, and $\epsilon_{\pm}(T)=\frac12(1\pm(-1)^{\bar(T)})$, which is either $1$ or $0$, 
while
\begin{equation}\label{eqn-Twgts-so2n}
\begin{array}{|l|l|l|l|}
\hline
T_{ij}&\wgt(T_{ij})&\cr
\hline
k&x_k+a_{2k-1+n-j-i+\delta_{ik}}&a_m=0\hbox{~~for~~}m\leq0\cr
 \ov{k}&\ov{x}_k+a_{2k+n-j-i}&\cr		
\hline	
\end{array}
\end{equation}																														
\end{Theorem}

\noindent{\bf Proof}:
The $\ell(\lambda)<n$ result (\ref{eqn-wgt-tab-o2n}) is clear since $so_\lambda(\x,\ovx|\a)=o_\lambda(\x,\ovx|\a)$
in this case, allowing us to use the result from Theorem~\ref{thm-Twgts}. For $\ell(\lambda)=n$
a factor $2^{\zeta(T)}$ with $\zeta>0$ is associated with the doubling that necessarily occurs following the action 
of $\psi$ in standardising tableaux that would otherwise have $\ov{k}$ appearing immediately above $k$ in the
first column for some $\zeta(T)$ values of $k$. These tableaux can only arise through the involvement of one or more type (iii) paths
in the corresponding non-intersecting $n$-tuple. They therefore contribute to $o_\lambda(\x,\ovx|\a)$ but not to $o'_\lambda(\x,\ovx|\a)$. 
The contribution to $so_{\lambda_\pm}(\x,\ovx|\a)=\frac12(o_\lambda(\x,\ovx|\a)\pm o'_\lambda(\x,\ovx|\a))$ must then be weighted 
by $2^{\zeta(T)-1}$, as required in (\ref{eqn-wgt-tab-so2n}). There just remain the terms with $\zeta(T)=0$ that contribute 
to both $o_\lambda(\x,\ovx|\a)$ and $o'_\lambda(\x,\ovx|\a)$, but to the latter with a sign factor $(-1)^{\bar(T)}$. It follows that
this contributes to $so_{\lambda_{\pm}}(\x,\ovx|\a)$ according as $\frac12(1\pm(-1)^{\bar(T)})$ is $1$ or $0$, as required in (\ref{eqn-wgt-tab-so2n}).
The proof is then completed by noting that the weights applying to edges and thus tableau entries coincide with those in the $o(2n)$ case.
\qed

By way of example, consider the case $n=2$ and $\lambda=(2,2)$ for which we
have the following contributions to the factorial irreducible characters of $SO(4)$ of highest weights $(2,2)$ and $(2,-2)$.

First, the case $so_{(2,2)_+}(\x,\ovx|\a)$
$$
\begin{array}{ccc}
\vcenter{\hbox{
\begin{tikzpicture}[x={(0in,-0.4in)},y={(0.5in,0in)}] 
\foreach \j in {2,...,4} \draw(2,\j)node{$\bullet$};
\foreach \j in {2,...,4} \draw(3,\j)node{$\bullet$};
\foreach \j in {1,...,4} \draw(4,\j)node{$\bullet$};
\draw(2-0.3,2-0.3)node{$P_1$};
\draw(4-0.3,1-0.3)node{$P_2$};
\draw(4+0.3,3+0.3)node{$Q_2$};
\draw(4+0.3,4+0.3)node{$Q_1$};
\draw[-](4,1)to(1.0,4.0); \draw(1.0,4.0+0.3)node{$a_{0}$};
\draw[-](4,2)to(1.5,4.5); \draw(1.5,4.5+0.3)node{$a_1$};
\draw[-](4,3)to(2.5,4.5); \draw(2.5,4.5+0.3)node{$a_2$};
\draw[-](4,4)to(3.5,4.5); \draw(3.5,4.5+0.3)node{$a_3$};
\draw[draw=blue,ultra thick] (2,2)to(2,3); \draw(2-0.3,3-0.4)node{$x_1\!\!+\!\!{a_{0}}$};
\draw[draw=blue,ultra thick] (2,3)to(2,4); \draw(2-0.3,4-0.4)node{$x_1\!\!+\!\!{a_{1}}$};
\draw[draw=blue,ultra thick] (2,4)to(4,4);
\draw[draw=red,ultra thick] (4,1)to(4,2); \draw(4-0.3,2-0.4)node{$x_2\!\!+\!\!{a_1}$};
\draw[draw=red,ultra thick] (4,2)to(4,3); \draw(4-0.3,3-0.4)node{$x_2\!\!+\!\!{a_2}$};
\end{tikzpicture}
}}
&
\vcenter{\hbox{
\begin{tikzpicture}[x={(0in,-0.4in)},y={(0.5in,0in)}] 
\foreach \j in {2,...,4} \draw(2,\j)node{$\bullet$};
\foreach \j in {2,...,4} \draw(3,\j)node{$\bullet$};
\foreach \j in {1,...,4} \draw(4,\j)node{$\bullet$};
\draw[-](4,1)to(1.0,4.0); \draw(1.0,4.0+0.3)node{$a_{0}$};
\draw[-](4,2)to(1.5,4.5); \draw(1.5,4.5+0.3)node{$a_1$};
\draw[-](4,3)to(2.5,4.5); \draw(2.5,4.5+0.3)node{$a_2$};
\draw[-](4,4)to(3.5,4.5); \draw(3.5,4.5+0.3)node{$a_3$};
\draw[draw=blue,ultra thick] (2,2)to(2,3); \draw(2-0.3,3-0.5)node{$x_1\!\!+\!\!{a_0}$};
\draw[draw=blue,ultra thick] (2,3)to(3,3); 
\draw[draw=blue,ultra thick] (3,3)to(3,4); \draw(3-0.3,4-0.5)node{$\ov{x}_2\!\!+\!\!a_2$};
\draw[draw=blue,ultra thick] (3,4)to(4,4); 
\draw[draw=red,ultra thick] (4,1)to(4,2); \draw(4-0.3,2-0.5)node{$x_2\!\!+\!\!{a_1}$};
\draw[draw=red,ultra thick] (4,2)to(4,3); \draw(4-0.3,3-0.5)node{${x}_2\!\!+\!\!{a_2}$};
\foreach \j in {2,...,4} \draw(2,\j)node{$\bullet$};
\foreach \j in {2,...,4} \draw(3,\j)node{$\bullet$};
\foreach \j in {1,...,4} \draw(4,\j)node{$\bullet$};
\draw(2-0.3,2-0.3)node{$P_1$};
\draw(4-0.3,1-0.3)node{$P_2$};
\draw(4+0.3,3+0.3)node{$Q_2$};
\draw(4+0.3,4+0.3)node{$Q_1$};
\end{tikzpicture}
}}
&
\cr\cr
\wideYT{0.2in}{0.6in}{}{
 {x_1\!+\!{a_0},x_1\!+\!{a_1}},
 {x_2\!+\!{a_1},x_2\!+\!{a_2}},
} 
&
\psi_2:
\wideYT{0.2in}{0.6in}{}{
 {x_1\!+\!{a_0},\ov{x}_2\!+\!a_2},
 {x_2\!+\!{a_1},{x}_2\!+\!{a_2}},
} \mapsto
\wideYT{0.2in}{0.6in}{}{
 {x_1\!+\!{a_0},{x}_2\!+\!a_2},
 {x_2\!+\!{a_1},\ov{x}_2\!+\!{a_2}},
} 
\end{array}
$$

$$
\begin{array}{cc}
\vcenter{\hbox{
\begin{tikzpicture}[x={(0in,-0.4in)},y={(0.5in,0in)}] 
\foreach \j in {2,...,4} \draw(2,\j)node{$\bullet$};
\foreach \j in {2,...,4} \draw(3,\j)node{$\bullet$};
\foreach \j in {1,...,4} \draw(4,\j)node{$\bullet$};
\draw[-](4,1)to(1.0,4.0); \draw(1.0,4.0+0.3)node{$a_{0}$};
\draw[-](4,2)to(1.5,4.5); \draw(1.5,4.5+0.3)node{$a_1$};
\draw[-](4,3)to(2.5,4.5); \draw(2.5,4.5+0.3)node{$a_2$};
\draw[-](4,4)to(3.5,4.5); \draw(3.5,4.5+0.3)node{$a_3$};
\draw[draw=blue,ultra thick] (2,2)to(2,3); \draw(2-0.3,3-0.5)node{$\ov{x}_1\!\!+\!\!a_{0}$};
\draw[draw=blue,ultra thick] (2,3)to(3,3);
\draw[draw=blue,ultra thick] (3,3)to(3,4); \draw(3-0.3,4-0.5)node{$x_2\!\!+\!\!a_2$};
\draw[draw=blue,ultra thick] (3,4)to(4,4); 
\draw[draw=red,ultra thick] (4,1)to(4,2); \draw(4-0.3,2-0.5)node{$\ov{x}_2\!\!+\!\!a_1$};
\draw[draw=red,ultra thick] (4,2)to(4,3); \draw(4-0.3,3-0.5)node{$\ov{x}_2\!\!+\!\!{a_2}$};
\foreach \j in {2,...,4} \draw(2,\j)node{$\bullet$};
\foreach \j in {2,...,4} \draw(3,\j)node{$\bullet$};
\foreach \j in {1,...,4} \draw(4,\j)node{$\bullet$};
\draw(2-0.3,2-0.3)node{$P_1$};
\draw(4-0.3,1-0.3)node{$P_2$};
\draw(4+0.3,3+0.3)node{$Q_2$};
\draw(4+0.3,4+0.3)node{$Q_1$};
\end{tikzpicture}
}}
&
\vcenter{\hbox{
\begin{tikzpicture}[x={(0in,-0.4in)},y={(0.5in,0in)}] 
\foreach \j in {2,...,4} \draw(2,\j)node{$\bullet$};
\foreach \j in {2,...,4} \draw(3,\j)node{$\bullet$};
\foreach \j in {1,...,4} \draw(4,\j)node{$\bullet$};
\draw[-](4,1)to(1.0,4.0); \draw(1.0,4.0+0.3)node{$a_{0}$};
\draw[-](4,2)to(1.5,4.5); \draw(1.5,4.5+0.3)node{$a_1$};
\draw[-](4,3)to(2.5,4.5); \draw(2.5,4.5+0.3)node{$a_2$};
\draw[-](4,4)to(3.5,4.5); \draw(3.5,4.5+0.3)node{$a_3$};
\draw[draw=blue,ultra thick] (2,2)to(2,3); \draw(2-0.3,3-0.5)node{$\ov{x}_1\!\!+\!\!a_0$};
\draw[draw=blue,ultra thick] (2,3)to(2,4); \draw(2-0.3,4-0.5)node{$\ov{x}_1\!\!+\!\!a_1$};
\draw[draw=blue,ultra thick] (2,4)to(4,4); 
\draw[draw=red,ultra thick] (4,1)to(4,2); \draw(4-0.3,2-0.5)node{$\ov{x}_2\!\!+\!\!{a_1}$};
\draw[draw=red,ultra thick] (4,2)to(4,3); \draw(4-0.3,3-0.5)node{$\ov{x}_2\!\!+\!\!{a_2}$};
\foreach \j in {2,...,4} \draw(2,\j)node{$\bullet$};
\foreach \j in {2,...,4} \draw(3,\j)node{$\bullet$};
\foreach \j in {1,...,4} \draw(4,\j)node{$\bullet$};
\draw(2-0.3,2-0.3)node{$P_1$};
\draw(4-0.3,1-0.3)node{$P_2$};
\draw(4+0.3,3+0.3)node{$Q_2$};
\draw(4+0.3,4+0.3)node{$Q_1$};
\end{tikzpicture}
}}
\cr\cr
\wideYT{0.2in}{0.6in}{}{
 {\ov{x}_1\!+\!a_{0},x_2\!+\!a_2},
 {\ov{x}_2\!+\!a_1,\ov{x}_2\!+\!{a_2}},
} 
&
\wideYT{0.2in}{0.6in}{}{
 {\ov{x}_1\!+\!{a_0},\ov{x}_1\!+\!a_1},
 {\ov{x}_2\!+\!{a_1},\ov{x}_2\!+\!{a_2}},
} 
\end{array}
$$

$$
\begin{array}{c}
\vcenter{\hbox{
\begin{tikzpicture}[x={(0in,-0.4in)},y={(0.5in,0in)}] 
\foreach \j in {2,...,4} \draw(2,\j)node{$\bullet$};
\foreach \j in {2,...,4} \draw(3,\j)node{$\bullet$};
\foreach \j in {1,...,4} \draw(4,\j)node{$\bullet$};
\draw[-](4,1)to(1.0,4.0); \draw(1.0,4.0+0.3)node{$a_{0}$};
\draw[-](4,2)to(1.5,4.5); \draw(1.5,4.5+0.3)node{$a_1$};
\draw[-](4,3)to(2.5,4.5); \draw(2.5,4.5+0.3)node{$a_2$};
\draw[-](4,4)to(3.5,4.5); \draw(3.5,4.5+0.3)node{$a_3$};
\draw[draw=blue,ultra thick] (2,2)to(3,2); 
\draw[draw=blue,ultra thick] (3,2)to(3,3); \draw(3-0.3,3-0.5)node{$x_2\!\!+\!\!a_{1}$};
\draw[draw=blue,ultra thick] (3,3)to(3,4); \draw(3-0.3,4-0.5)node{$x_2\!\!+\!\!a_2$};
\draw[draw=blue,ultra thick] (3,4)to(4,4); 
\draw[draw=red,ultra thick] (4,1)to(4,2); \draw(4-0.3,2-0.5)node{$\ov{x}_2\!\!+\!\!{a_1}$};
\draw[draw=red,ultra thick] (4,2)to(4,3); \draw(4-0.3,3-0.5)node{$\ov{x}_2\!\!+\!\!{a_2}$};
\foreach \j in {2,...,4} \draw(2,\j)node{$\bullet$};
\foreach \j in {2,...,4} \draw(3,\j)node{$\bullet$};
\foreach \j in {1,...,4} \draw(4,\j)node{$\bullet$};
\draw(2-0.3,2-0.3)node{$P_1$};
\draw(4-0.3,1-0.3)node{$P_2$};
\draw(4+0.3,3+0.3)node{$Q_2$};
\draw(4+0.3,4+0.3)node{$Q_1$};
\end{tikzpicture}
}}
\cr\cr
\wideYT{0.2in}{0.6in}{}{
 {x_2\!+\!{a_1},x_2\!+\!a_2},
 {\ov{x}_2\!+\!{a_1},\ov{x}_2\!+\!{a_2}},
} 
\end{array}
$$
The factorial character $so_{(2,2)_+}(\x,\ovx|\a)$ is then obtained by summing these five contributions and
setting $a_0=0$.

Then the case $so_{(2,2)_-}(\x,\ovx|\a)$
$$
\begin{array}{ccc}
\vcenter{\hbox{
\begin{tikzpicture}[x={(0in,-0.4in)},y={(0.5in,0in)}] 
\foreach \j in {2,...,4} \draw(2,\j)node{$\bullet$};
\foreach \j in {2,...,4} \draw(3,\j)node{$\bullet$};
\foreach \j in {1,...,4} \draw(4,\j)node{$\bullet$};
\draw[-](4,1)to(1.0,4.0); \draw(1.0,4.0+0.3)node{$a_{0}$};
\draw[-](4,2)to(1.5,4.5); \draw(1.5,4.5+0.3)node{$a_1$};
\draw[-](4,3)to(2.5,4.5); \draw(2.5,4.5+0.3)node{$a_2$};
\draw[-](4,4)to(3.5,4.5); \draw(3.5,4.5+0.3)node{$a_3$};
\draw[draw=blue,ultra thick] (2,2)to(2,3); \draw(2-0.3,3-0.5)node{$x_1\!\!+\!\!{a_0}$};
\draw[draw=blue,ultra thick] (2,3)to(2,4); \draw(2-0.3,4-0.5)node{$x_1\!\!+\!\!{a_1}$};
\draw[draw=blue,ultra thick] (2,4)to(4,4);
\draw[draw=red,ultra thick] (4,1)to(4,2); \draw(4-0.3,2-0.5)node{$\ov{x}_2\!\!+\!\!{a_1}$};
\draw[draw=red,ultra thick] (4,2)to(4,3); \draw(4-0.3,3-0.5)node{$\ov{x}_2\!\!+\!\!{a_2}$};
\foreach \j in {2,...,4} \draw(2,\j)node{$\bullet$};
\foreach \j in {2,...,4} \draw(3,\j)node{$\bullet$};
\foreach \j in {1,...,4} \draw(4,\j)node{$\bullet$};
\draw(2-0.3,2-0.3)node{$P_1$};
\draw(4-0.3,1-0.3)node{$P_2$};
\draw(4+0.3,3+0.3)node{$Q_2$};
\draw(4+0.3,4+0.3)node{$Q_1$};
\end{tikzpicture}
}}
&
\vcenter{\hbox{
\begin{tikzpicture}[x={(0in,-0.4in)},y={(0.5in,0in)}] 
\foreach \j in {2,...,4} \draw(2,\j)node{$\bullet$};
\foreach \j in {2,...,4} \draw(3,\j)node{$\bullet$};
\foreach \j in {1,...,4} \draw(4,\j)node{$\bullet$};
\draw[-](4,1)to(1.0,4.0); \draw(1.0,4.0+0.3)node{$a_{0}$};
\draw[-](4,2)to(1.5,4.5); \draw(1.5,4.5+0.3)node{$a_1$};
\draw[-](4,3)to(2.5,4.5); \draw(2.5,4.5+0.3)node{$a_2$};
\draw[-](4,4)to(3.5,4.5); \draw(3.5,4.5+0.3)node{$a_3$};
\draw[draw=blue,ultra thick] (2,2)to(2,3); \draw(2-0.3,3-0.5)node{$x_1\!\!+\!\!{a_0}$};
\draw[draw=blue,ultra thick] (2,3)to(3,3); 
\draw[draw=blue,ultra thick] (3,3)to(3,4); \draw(3-0.3,4-0.5)node{$x_2\!\!+\!\!a_2$};
\draw[draw=blue,ultra thick] (3,4)to(4,4); 
\draw[draw=red,ultra thick] (4,1)to(4,2); \draw(4-0.3,2-0.5)node{$\ov{x}_2\!\!+\!\!a_1$};
\draw[draw=red,ultra thick] (4,2)to(4,3); \draw(4-0.3,3-0.5)node{$\ov{x}_2\!\!+\!\!{a_2}$};
\foreach \j in {2,...,4} \draw(2,\j)node{$\bullet$};
\foreach \j in {2,...,4} \draw(3,\j)node{$\bullet$};
\foreach \j in {1,...,4} \draw(4,\j)node{$\bullet$};
\draw(2-0.3,2-0.3)node{$P_1$};
\draw(4-0.3,1-0.3)node{$P_2$};
\draw(4+0.3,3+0.3)node{$Q_2$};
\draw(4+0.3,4+0.3)node{$Q_1$};
\end{tikzpicture}
}}
&
\cr\cr
\wideYT{0.2in}{0.6in}{}{
 {x_1\!+\!{a_0},{x}_1\!+\!{a_1}},
 {\ov{x}_2\!+\!{a_1},\ov{x}_2\!+\!{a_2}},
} 
&
\wideYT{0.2in}{0.6in}{}{
 {x_1\!+\!{a_0},x_2\!+\!a_2},
 {\ov{x}_2\!+\!{a_1},\ov{x}_2\!+\!{a_2}},
} 
\end{array}
$$

$$
\begin{array}{cc}
\vcenter{\hbox{
\begin{tikzpicture}[x={(0in,-0.4in)},y={(0.5in,0in)}] 
\foreach \j in {2,...,4} \draw(2,\j)node{$\bullet$};
\foreach \j in {2,...,4} \draw(3,\j)node{$\bullet$};
\foreach \j in {1,...,4} \draw(4,\j)node{$\bullet$};
\draw[-](4,1)to(1.0,4.0); \draw(1.0,4.0+0.3)node{$a_{0}$};
\draw[-](4,2)to(1.5,4.5); \draw(1.5,4.5+0.3)node{$a_1$};
\draw[-](4,3)to(2.5,4.5); \draw(2.5,4.5+0.3)node{$a_2$};
\draw[-](4,4)to(3.5,4.5); \draw(3.5,4.5+0.3)node{$a_3$};
\draw[draw=blue,ultra thick] (2,2)to(2,3); \draw(2-0.3,3-0.5)node{$\ov{x}_1\!\!+\!\!a_{0}$};
\draw[draw=blue,ultra thick] (2,3)to(3,3);
\draw[draw=blue,ultra thick] (3,3)to(3,4); \draw(3-0.3,4-0.5)node{$\ov{x}_2\!\!+\!\!a_2$};
\draw[draw=blue,ultra thick] (3,4)to(4,4); 
\draw[draw=red,ultra thick] (4,1)to(4,2); \draw(4-0.3,2-0.5)node{$x_2\!\!+\!\!{a_1}$};
\draw[draw=red,ultra thick] (4,2)to(4,3); \draw(4-0.3,3-0.5)node{${x}_2\!\!+\!\!{a_2}$};
\foreach \j in {2,...,4} \draw(2,\j)node{$\bullet$};
\foreach \j in {2,...,4} \draw(3,\j)node{$\bullet$};
\foreach \j in {1,...,4} \draw(4,\j)node{$\bullet$};
\draw(2-0.3,2-0.3)node{$P_1$};
\draw(4-0.3,1-0.3)node{$P_2$};
\draw(4+0.3,3+0.3)node{$Q_2$};
\draw(4+0.3,4+0.3)node{$Q_1$};
\end{tikzpicture}
}}
&
\vcenter{\hbox{
\begin{tikzpicture}[x={(0in,-0.4in)},y={(0.5in,0in)}] 
\foreach \j in {2,...,4} \draw(2,\j)node{$\bullet$};
\foreach \j in {2,...,4} \draw(3,\j)node{$\bullet$};
\foreach \j in {1,...,4} \draw(4,\j)node{$\bullet$};
\draw[-](4,1)to(1.0,4.0); \draw(1.0,4.0+0.3)node{$a_{0}$};
\draw[-](4,2)to(1.5,4.5); \draw(1.5,4.5+0.3)node{$a_1$};
\draw[-](4,3)to(2.5,4.5); \draw(2.5,4.5+0.3)node{$a_2$};
\draw[-](4,4)to(3.5,4.5); \draw(3.5,4.5+0.3)node{$a_3$};
\draw[draw=blue,ultra thick] (2,2)to(2,3); \draw(2-0.3,3-0.5)node{$\ov{x}_1\!\!+\!\!a_0$};
\draw[draw=blue,ultra thick] (2,3)to(2,4); \draw(2-0.3,4-0.5)node{$\ov{x}_1\!\!+\!\!a_1$};
\draw[draw=blue,ultra thick] (2,4)to(4,4); 
\draw[draw=red,ultra thick] (4,1)to(4,2); \draw(4-0.3,2-0.5)node{$x_2\!\!+\!\!{a_1}$};
\draw[draw=red,ultra thick] (4,2)to(4,3); \draw(4-0.3,3-0.5)node{$x_2\!\!+\!\!{a_2}$};
\foreach \j in {2,...,4} \draw(2,\j)node{$\bullet$};
\foreach \j in {2,...,4} \draw(3,\j)node{$\bullet$};
\foreach \j in {1,...,4} \draw(4,\j)node{$\bullet$};
\draw(2-0.3,2-0.3)node{$P_1$};
\draw(4-0.3,1-0.3)node{$P_2$};
\draw(4+0.3,3+0.3)node{$Q_2$};
\draw(4+0.3,4+0.3)node{$Q_1$};
\end{tikzpicture}
}}
\cr\cr
\psi_2:\wideYT{0.2in}{0.6in}{}{
 {\ov{x}_1\!+\!a_{0},\ov{x}_2\!+\!a_2},
 {x_2\!+\!{a_1},{x}_2\!+\!{a_2}},
} \mapsto
\wideYT{0.2in}{0.6in}{}{
 {\ov{x}_1\!+\!a_{0},{x}_2\!+\!a_2},
 {x_2\!+\!{a_1},\ov{x}_2\!+\!{a_2}},
}
&
\wideYT{0.2in}{0.6in}{}{
 {\ov{x}_1\!+\!a_{0},\ov{x}_1\!+\!a_1},
 {x_2\!+\!{a_1},x_2\!+\!{a_2}},
} 
\end{array}
$$

$$
\begin{array}{c}
\vcenter{\hbox{
\begin{tikzpicture}[x={(0in,-0.4in)},y={(0.5in,0in)}] 
\foreach \j in {2,...,4} \draw(2,\j)node{$\bullet$};
\foreach \j in {2,...,4} \draw(3,\j)node{$\bullet$};
\foreach \j in {1,...,4} \draw(4,\j)node{$\bullet$};
\draw[-](4,1)to(1.0,4.0); \draw(1.0,4.0+0.3)node{$a_{0}$};
\draw[-](4,2)to(1.5,4.5); \draw(1.5,4.5+0.3)node{$a_1$};
\draw[-](4,3)to(2.5,4.5); \draw(2.5,4.5+0.3)node{$a_2$};
\draw[-](4,4)to(3.5,4.5); \draw(3.5,4.5+0.3)node{$a_3$};
\draw[draw=blue,ultra thick] (2,2)to(3,2); 
\draw[draw=blue,ultra thick] (3,2)to(3,3); \draw(3-0.3,3-0.5)node{$\ov{x}_2\!\!+\!\!a_{1}$};
\draw[draw=blue,ultra thick] (3,3)to(3,4); \draw(3-0.3,4-0.5)node{$\ov{x}_2\!\!+\!\!a_2$};
\draw[draw=blue,ultra thick] (3,4)to(4,4); 
\draw[draw=red,ultra thick] (4,1)to(4,2); \draw(4-0.3,2-0.5)node{${x}_2\!\!+\!\!{a_1}$};
\draw[draw=red,ultra thick] (4,2)to(4,3); \draw(4-0.3,3-0.5)node{${x}_2\!\!+\!\!{a_2}$};
\foreach \j in {2,...,4} \draw(2,\j)node{$\bullet$};
\foreach \j in {2,...,4} \draw(3,\j)node{$\bullet$};
\foreach \j in {1,...,4} \draw(4,\j)node{$\bullet$};
\draw(2-0.3,2-0.3)node{$P_1$};
\draw(4-0.3,1-0.3)node{$P_2$};
\draw(4+0.3,3+0.3)node{$Q_2$};
\draw(4+0.3,4+0.3)node{$Q_1$};
\end{tikzpicture}
}}
\cr\cr
\psi_2:\wideYT{0.2in}{0.6in}{}{
 {\ov{x}_2\!+\!{a_1},\ov{x}_2\!+\!a_2},
 {{x}_2\!+\!{a_1},{x}_2\!+\!{a_2}},
} \mapsto
\wideYT{0.2in}{0.6in}{}{
 {x_2\!+\!{a_1},x_2\!+\!a_2},
 {\ov{x}_2\!+\!{a_1},\ov{x}_2\!+\!{a_2}},
}
\end{array}
$$
The factorial character $so_{(2,2)_-}(\x,\ovx|\a)$ is then obtained by summing these five contributions and
setting $a_0=0$.

It should be noted that we have only indicated the action of $\psi$ where the original tableau, obtained directly
from the $2$-tuple of non-intersecting lattice paths, is not in ${\cal T}_{(2,2)}^{\eo}$. In this example,
it is only $\psi_2$ that need be applied and it only has the effect of interchanging entries in vertical pairs. 
The contributions to $so_{(2,2)_+}(\x,\ovx|\a)$ and $so_{(2,2)_-}(\x,\ovx|\a)$ arise only from those $2$-tuples
containing an even and an odd number, respectively, of type (ii) paths. In this case the action of $\psi_2$
renders the fifth tableau contributing to $so_{(2,2)_-}(\x,\ovx|\a)$ identical to the fifth tableau contributing 
to $so_{(2,2)_+}(\x,\ovx|\a)$, thereby giving rise to a multiplicity of $2$ in the even orthogonal character $o_{(2,2)}(\x,\ovx|\a)$. 


\section{Closing remarks}
\label{sec-remarks}

Our results are based heavily on the Definition~\ref{def-hma} of $h^{\g}_m(\z|\a)$ in terms of 
generating functions that are manifestly Weyl symmetric. However, these symmetries are not
always evident in our main results. For example, in Theorem~\ref{thm-fJT} we have chosen 
to express the result in terms of a particular choice of flag, namely one for which
$\x^{(i)}=(x_i,x_{i+1},\ldots,x_n)$ with $\x^{(1)}\supset \x^{(2)} \supset \cdots \supset \x^{(n)}$
and $\x^{(i)}\backslash\x^{(i+1)}=x_i$ for $i=1,2,\ldots,n-1$ and $\x^{(n)}=x_n$. 
However, the overall symmetry with respect to permutations of $\x$ means that the results 
are independent of this particular choice of flag. 
In particular one might equally well define $\x_{(i)}=(x_1,x_2,\ldots,x_{i})$ and adopt a flag
$\x_{(1)}\subset \x_{(2)} \subset \cdots \subset \x_{(n)}$ with $\x_{(1)}=x_1$ and 
$\x_{(i)}\backslash/\x_{(i-1)}=x_i$ for $i=2,3,\ldots,n$.

This freedom of choice is particularly important when it comes to labelling rows in our
lattice path models and ordering entries in our tableau models. In general we have opted for models leading to rather
well known sets of tableaux such as those based on the order $1<2<\cdots<n$ familiar in the definition 
of semistandard tableaux. When it came to the symplectic case, for which $h^{\sp}_m(\x,\ovx|\a)$
is symmetric under the interchange of any $x_i$ and $\ov{x}_j$, we have chosen to 
interleave barred and unbarred elements with an order $1<\ov1<2<\ov2<\cdots<n<\ov{n}$ 
as used for example by Koike and Terada~\cite{KT90} which is 
obviously very little different from the use of $\ov1<1<\ov2<2<\cdots<\ov{n}<n$ in~\cite{Kin76}.
However, we might have tried to adopt the order $\ov{n}<\cdots<\ov2<\ov1<1<2<\cdots< n$, 
as used by De Concini~\cite{DeC79}, or its variant $1<2<\cdots n<\ov{n}<\cdots<\ov2<\ov1$ as used         
in crystal graphs by Kashiwara and Nakashima~\cite{KN94}. There would then remain the task of devising suitable
weightings of the tableaux based on these orders so as to yield factorial symplectic characters.

Similarly, in the odd orthogonal case we have deliberately chosen the order $1<\ov1<2<\ov2<\cdots<n<\ov{n}<0$ so as
to recreate the tableaux of Sundaram~\cite{Sun90}, albeit with her symbol $\infty$ replaced by $0$.
Alternatively, we could have adopted the alphabet $\sharp_1<1<\ov1<\sharp_2<2<\ov2<\cdots<\sharp_n<n<\ov{n}$
introduced by Koike and Terada~\cite{KT90} in specifying their $so(2n+1)$ tableaux. In this case, the factorial weights 
attached to entries $k$, $\ov{k}$ and $\sharp_i$ at position $(i,j)$, $(i,j)$ and $(i,1)$ are 
$(x_k+a_{3k-1-2n+j-i})$, $(\ov{x}_k+a_{3k-2n+j-i})$ and $(1-a_{2i-2n})$, respectively, with $a_m=0$ for $m\leq0$. 
Koike and Terada also introduced both $\sharp_i$ and $\flat_i$ in the same way in specifying their $so(2n)$ tableaux. 
The weightings of such tableaux that are necessary to give factorial even orthogonal characters remain to be determined.

While all this supports the case for defining factorial characters as in Section~\ref{sec-fchar}, in order to go
further and arrive at appropriate Tokuyama-type identities for such characters analagous to those obtained already for the factorial
Schur function case~\cite{BMN14,HK15b}, it is necessary to define factorial $Q$-functions $Q^{\g}_\lambda(\w;\z|\a)$ for each $\g$
and strict partition $\lambda$ of length $\ell(\lambda)\leq n$. This will be the subject of further paper in which each 
$Q^{\g}_\lambda(\w;\z|\a)$ is defined in terms of determinants whose elements are simple multiples of supersymmetric factorial 
$q$-functions $q^{\g}_m(\w;\z|\a)$, which are themselves defined by means of generating functions. 
These algebraic definitions of $Q^{\g}_\lambda(\w;\z|\a)$ may then be converted into combinatorial expressions by following the same procedure 
as that used here for $\g_\lambda(\z|\a)$and $h_m(\z|\a)$, that is proceeding by way of a lattice path model to a model
involving this time sets, ${\cal P}^{\g}_\lambda$, of primed shifted tableaux $P$ of shifted shape $\lambda$.

\bigskip


\noindent{\bf Acknowledgements}
\label{sec:ack}
This work was supported by the Canadian Tri-Council Research
Support Fund. The first author (AMH) acknowledges the support of a
Discovery Grant from the Natural Sciences and Engineering Research Council of
Canada (NSERC). The second author (RCK) is grateful for the hospitality extended to him
by AMH and her colleagues at Wilfrid Laurier University, by Professors Bill Chen and Arthur Yang at Nankai University and 
by Professor Itaru Terada at the University of Tokyo, and for the financial support making visits to each of these
universities possible.

\vfill
\end{document}